\documentclass[11pt]{article}
\usepackage{amssymb,amsmath,latexsym, verbatim, enumerate,color}
\allowdisplaybreaks

\begin{document}

\begin{doublespace}

\newtheorem{thm}{Theorem}[section]
\newtheorem{lemma}[thm]{Lemma}
\newtheorem{cond}[thm]{Condition}
\newtheorem{defn}[thm]{Definition}
\newtheorem{prop}[thm]{Proposition}
\newtheorem{corollary}[thm]{Corollary}
\newtheorem{remark}[thm]{Remark}
\newtheorem{example}[thm]{Example}
\newtheorem{conj}[thm]{Conjecture}
\numberwithin{equation}{section}
\def\ee{\varepsilon}
\def\qed{{\hfill $\Box$ \bigskip}}
\def\NN{{\cal N}}
\def\AA{{\cal A}}
\def\MM{{\cal M}}
\def\BB{{\cal B}}
\def\CC{{\cal C}}
\def\LL{{\cal L}}
\def\DD{{\cal D}}
\def\FF{{\cal F}}
\def\EE{{\cal E}}
\def\QQ{{\cal Q}}
\def\RR{{\mathbb R}}
\def\R{{\mathbb R}}
\def\L{{\bf L}}
\def\K{{\bf K}}
\def\S{{\bf S}}
\def\A{{\bf A}}
\def\E{{\mathbb E}}
\def\F{{\bf F}}
\def\P{{\mathbb P}}
\def\N{{\mathbb N}}
\def\eps{\varepsilon}
\def\wh{\widehat}
\def\wt{\widetilde}
\def\pf{\noindent{\bf Proof.} }
\def\beq{\begin{equation}}
\def\eeq{\end{equation}}
\def\lam{\lambda}
\def\H{\mathcal{H}}
\def\nn{\nonumber}
\def\L{\mathcal{L}}

\def\HH{{\mathbb H}}

\newcommand{\Per}{\mathrm{Per}}
\newcommand{\norm}[1]{{\lVert #1 \rVert}}

\title{\Large \bf  Spectral heat content for time-changed killed Brownian motions}
\author{ \bf  Kei Kobayashi and Hyunchul Park }

\date{ \today}
\maketitle

\begin{abstract}
The spectral heat content is investigated for time-changed killed Brownian motions on $C^{1,1}$ open sets, where the time change is given by either a subordinator or an inverse subordinator, with the underlying Laplace exponent being regularly varying at $\infty$ with index $\beta\in (0,1)$. 
In the case of inverse subordinators, the asymptotic limit of the spectral heat content in small time is shown to involve a probabilistic term depending only on $\beta\in (0,1)$. In contrast, in the case of subordinators, this universality holds only when $\beta\in (\frac12, 1)$. 
\end{abstract}

\section{Introduction}\label{section:intro}
Consider the following heat equation with a Dirichlet boundary condition and an initial condition:
$$
\begin{cases}
\partial_t u(x,t)=\Delta u(x,t), & x\in \Omega,\ t>0,\\
u(x,0)=1, &x \in \Omega,\\
u(y,t)=0, & y\in \partial \Omega,\ t>0. 
\end{cases}
$$
The Laplace operator $\Delta$ is the infinitesimal generator of a Brownian motion, and the spectral heat content for the Brownian motion is defined by
$$
Q_{\Omega}^{BM}(t):=\int_{\Omega}u(x,t)dx=\int_{\Omega}\P_{x}(\tau_{\Omega}^{BM}>t)dx,
$$
where $\tau_{\Omega}^{BM}$ is the first exit time of the Brownian motion from the domain $\Omega$.
Intuitively, the spectral heat content measures the total heat that remains in the domain $\Omega$ at time $t>0$. 
The spectral heat content for the Brownian motion has been studied intensively in the past. One interesting feature is that the spectral heat content for a smooth domain admits a complete asymptotic expansion whose coefficients provide geometric information about the domain. In particular, the first and second terms involve the Lebesgue measure $|\Omega|$ and the perimeter $|\partial \Omega|$, respectively (\cite[Theorem 6.2]{vandenBerg1989}), while the third term depends on the mean curvature of $\Omega$ (\cite[Theorem 1.1]{BergLeGall94}).

Replacing the above Laplace operator $\Delta$ with the infinitesimal generator of a general L\'evy process and changing the zero boundary condition into the zero exterior condition (i.e.,\ $u(y,t)=0$ for all $y\in \Omega^{c}$) leads to the definition of the spectral heat content for the L\'evy process. 
The spectral heat contents for killed subordinate Brownian motions, which form a subclass of killed L\'evy processes, have been recently studied in \cite{GPS19, Valverde2017}, while the regular heat contents for more general L\'evy processes have been investigated in \cite{CG17}.

When a L\'evy process is given by a subordinate Brownian motion, which is a Brownian motion composed with an independent subordinator, the corresponding spectral heat content is the spectral heat content for the \textit{killed subordinate Brownian motion}. 
When one reverses the order of killing and time change (subordination), one can naturally consider the spectral heat content for a \textit{subordinate killed Brownian motion}. This is the main object of interest in this paper.
In fact, the paper deals with not only time-changed Brownian motions whose time changes are \textit{subordinators} but also those whose time changes are \textit{inverse subordinators}. 
The latter processes are important both theoretically and practically since they naturally appear in the context of subdiffusions, where the particles spread at a slower pace than the regular Brownian particles. 
In the past few decades, such time-changed processes and their various extensions have been widely studied from different perspectives (see e.g.\ \cite{JinKobayashi,Magdziarz_spa,MagdziarzZorawik,MNV1,MS_1,NaneNi2018,HKU-book} and references therein).

Let us stress that
the spectral heat content for time-changed killed Brownian motions is a natural object to study, and it also provides useful information about the spectral heat content for killed time-changed Brownian motions. 
In fact, in \cite{ParkSong21} the asymptotic behavior of the spectral heat content for subordinate killed Brownian motions with respect to stable subordinators plays an essential role in establishing that for killed stable processes. We believe that results in this paper play a similar role when studying the asymptotic behavior of the spectral heat content for killed subordinate Brownian motions with respect to general subordinators.

This paper can be regarded as a natural continuation of \cite{ParkSong19}, but we emphasize that the results to be presented here go much beyond \cite{ParkSong19} in two aspects. 
First, whereas \cite{ParkSong19} investigated the asymptotic limit as $t\downarrow 0$ of the spectral heat content for subordinate killed Brownian motions 
when the underlying subordinator is \textit{stable}, this paper encompasses much more general time changes, including subordinators that are not necessarily stable and their inverses (or inverse subordinators). Note that the major difficulty of dealing with non-stable subordinators and their inverses lies in the fact that one can no longer rely on a scaling property (or self-similarity), which plays a crucial role in \cite{ParkSong19}. Consequently, the results to be established in this paper are based on various non-trivial asymptotic estimates that allow us to avoid the use of the scaling property. 

Our main assumption about the underlying subordinator, whether the time change is the subordinator itself or its inverse, is that the associated Laplace exponent $\phi$ is regularly varying at $\infty$ with index $\beta \in (0,1)$.
In terms of time changes given by subordinators, the three separate cases $\beta \in (\frac12,1)$, $\beta=\frac12$ and $\beta\in (0,\frac12)$ result in three statements of different nature about the spectral heat content (Theorems \ref{thm:stable2}, \ref{thm:withlowerorder} and \ref{thm:lessthan1/2}). 
We emphasize that the methods and tools for handling these three cases are significantly different. 
When $\beta\in (\frac12,1)$, in Lemma \ref{lemma:stable1}
we employ an \textit{approximate scaling property} of the subordinator, which was used in \cite[Lemmas 4.5 and 4.6]{GPS19}. When $\beta=\frac12$, in the setting where the subordinator $(D_t)$ is an independent sum of a $\frac12$-stable subordinator $(S_{t}^{(1/2)})$ and a certain subordinator representing perturbations, we utilize the key fact that $\E[\sqrt{D_{t}}1_{\{D_{t}\leq \delta\}}]$ is comparable to $\E[\sqrt{S_{t}^{(1/2)}}1_{\{S_{t}^{(1/2)}\leq \delta\}}]$ for any fixed constant $\delta>0$, which is derived in Lemma \ref{lemma:mixedconv}. 
On the other hand, the case when $\beta\in (0,\frac12)$ appeals to Proposition \ref{prop:1/2}, which concerns weak convergence of L\'evy measures. 

In contrast, when the time change is given by an inverse subordinator, we obtain a single statement for all values of $\beta\in(0,1)$ (Theorem \ref{thm:inverse_subordinator}). 
Here, the major technical result is Proposition \ref{prop:E_t}, which we prove using both de Bruijn's and Karamata's Tauberian theorems and their associated monotone density theorems. 
When the time change is an inverse stable subordinator, we even derive a complete asymptotic expansion of the spectral heat content with the coefficients explicitly written (Theorem \ref{thm:complete expansion});
this can be regarded as an analogue of a similar statement known for Brownian motions without random time changes. 
To the authors' knowledge, this is the first paper in the literature which discusses the spectral and regular heat contents for time-changed Brownian motions when the time changes are \textit{inverse} subordinators; thus, Theorems \ref{thm:inverse_subordinator} and  \ref{thm:complete expansion} provide a new perspective for the study of subdiffusions.

The second factor that makes this paper substantially different from \cite{ParkSong19} is that it reveals a \textit{universality} in the asymptotic limit of the spectral heat content for time-changed Brownian motions when the time changes are either subordinators $(D_t)$ with $\beta\in(\frac 12,1)$ or inverse subordinators $(E_t)$ with $\beta\in(0,1)$ (Theorems \ref{thm:stable2} and \ref{thm:inverse_subordinator}). 
Here, the \textit{universality} means that the asymptotic limit depends only on the index $\beta$ but not on the entire information about the subordinator $(D_{t})$ or the inverse subordinator $(E_{t})$.
In particular, the asymptotic limit only involves a probabilistic term expressed as the supremum of a one-dimensional Brownian motion observed up to $S_{1}^{(\beta)}$ or $E_{1}^{(\beta)}$, where $(S_{t}^{(\beta)})$ and $(E_{t}^{(\beta)})$ are a $\beta$-stable subordinator and its inverse, respectively, rather than the original (more general) time changes $(D_t)$ and $(E_t)$. 
In contrast, Theorem \ref{thm:lessthan1/2} shows that when the time changes are subordinators $(D_t)$ with $\beta\in (0,\frac12)$, the asymptotic limit depends on the underlying L\'evy measure. The threshold case $\beta=\frac 12$ is covered in Theorem \ref{thm:withlowerorder}; however, our result is limited to a special class of Laplace exponents and does not provide a complete picture.

The paper is organized as follows.
Section \ref{section:general setting} introduces the concepts of spectral and regular heat contents for time-changed Brownian motions and derives a theorem concerning those quantities in a general setting (Proposition \ref{prop:general_timechange}). Section \ref{section:SKBM} studies the spectral heat content for subordinate killed Brownian motions with the time changes being subordinators whose Laplace exponents are regularly varying at $\infty$ with index $\beta\in(0,1)$; 
we divide the section into three subsections to address the three separate cases $\beta \in (\frac12,1)$, $\beta=\frac12$, and $\beta\in (0,\frac12)$ (Theorems \ref{thm:stable2}, \ref{thm:withlowerorder}, and \ref{thm:lessthan1/2}). 
On the other hand, Section \ref{section:inverse subordinators} deals with time changes given by inverse subordinators and proves the universality result (Theorem \ref{thm:inverse_subordinator}). 
The latter section also derives a complete asymptotic expansion for the spectral heat content under the assumption that the time change is an inverse stable subordinator (Theorem \ref{thm:complete expansion}). 
Finally, Section \ref{section:examples} is devoted to applications of the above theorems, exhibiting some concrete examples of the asymptotic limits of spectral and regular heat contents.

\section{Spectral and regular heat contents for time-changed Brownian motions}\label{section:general setting}
This section studies the asymptotic behaviors as $t\downarrow 0$ of the spectral and regular heat contents for a general time-changed killed Brownian motion, where the time change is given by a stochastic process starting at 0 which may or may not have nondecreasing sample paths. 
The following notations are used throughout the paper. For a stochastic process $(X_t)_{t\ge 0}$ in $\R^d$ and an open set $\Omega$ in $\R^d$, the random time
\[
	\tau^X_\Omega=\inf\{ t>0: X_t\in \Omega^c\}
\]
is the first exit time of $(X_t)$ from $\Omega$. If $(X_t)$ is a Markov process with a family of probability measures $(\P^X_x)_{x\in \R^d}$ specifying its initial point, $(X^\Omega_t)_{t\ge 0}$ denotes the corresponding killed process defined by 
\[
	X^\Omega_t=
	\begin{cases} X_t \ \ &\textrm{if} \ t<\tau^X_\Omega,\\
			 \partial \ \ &\textrm{if} \ t\ge \tau^X_\Omega,
	\end{cases}
\]
where $\partial$ is a cemetery state. 
Consider a random time change given by an independent nondecreasing process  $(U_t)_{t\ge 0}$ in $\R_+:=[0,\infty)$ such that a.s., $U_0=0$ and $U_{t}> 0$ for all $t>0$.
Note that the time-changed processes $X\circ U:=(X_{U_t})_{t\ge 0}$ and $X^\Omega\circ U:= (X^\Omega_{U_t})_{t\ge 0}$ start at a point $x$ in $\R^d$ if and only if the outer process $(X_t)$ starts at $x$. 
In this paper, we set up the independent processes $(X_t)$ and $(U_t)$ on a product space with product probability measure $\P_x:=\P^X_x\times \P^U$ with 
obvious notations, and the corresponding expectations are denoted by $\E_x$, $\E^X_x$, and $\E^U$, respectively.
The spectral heat content of the time-changed process $(X_{U_t})_{t\ge 0}$ and that of the corresponding subordinate killed process $(X^\Omega_{U_t})_{t\ge 0}$ are respectively given by 
\begin{align}\label{0053}
	Q^{X\circ U}_\Omega(t):=\int_\Omega \P_x(\tau^{X\circ U}_\Omega>t)dx 
	\ \ \textrm{and} \ \ 
	\tilde{Q}^{X\circ U}_\Omega(t):=\int_\Omega \P_x(\tau^{X^\Omega\circ U}_\Omega>t)dx. 
\end{align}
The first quantity $Q^{X\circ U}_\Omega(t)$ concerns the killed subordinate process (where the subordinate process $(X_{U_t})_{t\ge 0}$ is killed upon exiting $\Omega$), while the second quantity $\tilde{Q}^{X\circ U}_\Omega(t)$ concerns the subordinate killed process (where the process $(X_t)$ is first killed upon exiting $\Omega$ and then subordinated to $(U_t)$). The two quantities satisfy the inequality $\tilde{Q}^{X\circ U}_\Omega(t)\le Q^{X\circ U}_\Omega(t)$ for any $t\ge 0$; see e.g.\ \cite{ParkSong19} for details. On the other hand, 
the regular heat content of $(X_{U_t})$ in $\Omega$ at time $t$ is given by
\[
	\HH^{X\circ U}_\Omega(t):=\int_\Omega\P_x(X_{U_t}\in \Omega)dx. 
\]
In particular, if $U_t=t$ for all $t\ge 0$ and $(X_t)$ is given by a Brownian motion $(W_t)$ with $\E_0[e^{-i\langle \xi, W_t\rangle}]=e^{-t|\xi|^2}$, then the above notions of the spectral and regular heat contents for the time-changed processes become those for $(W_t)$, and the following results hold 
for a bounded open interval $\Omega$ in $\R^1$ or a bounded 
connected
$C^{1,1}$ open set $\Omega$ in $\R^d$ with $d\ge 2$
(see \cite[Theorem 6.2]{vandenBerg1989}, \cite[Theorem 2]{vandenBerg2015}, and \cite[Theorem 2.4]{Miranda2007}):
\begin{align}
	\label{0031}\lim_{t\downarrow 0}\frac{|\Omega|-Q_{\Omega}^{W}(t)}{\sqrt{t}}&= \frac{2 |\partial \Omega|}{\sqrt{\pi}},\\
	\label{0032}\lim_{t\downarrow 0}\frac{|\Omega|-\HH^{W}_\Omega(t)}{\sqrt{t}}&= \frac{|\partial\Omega|}{\sqrt{\pi}},
\end{align}
where $|\Omega|$ and $|\partial\Omega|$ denote the $d$-dimensional Lebesgue measure of $\Omega$ and the $(d-1)$-dimensional Lebesgue measure of the boundary $\partial\Omega$, respectively. 
The reader is alerted that, throughout the paper, $(W_t)$ denotes a Brownian motion with $\E_0[e^{-i\langle \xi, W_t\rangle}]=e^{-t|\xi|^2}$, the governing equation of which is $\partial_t u(t,x) = \Delta u(t,x)$ rather than $\partial_t u(t,x) = \frac 12\Delta u(t,x)$.

The next lemma confirms that the above definitions of the spectral and regular heat contents are meaningful for a given time-changed Brownian motion $(W_{U_t})$. 
In the discussion below, we consider $(W_t)$ on the canonical space of continuous functions.

\begin{lemma}
Let $(W_t)$ be a Brownian motion in $\R^d$ independent of a process $(U_t)$ in $\R_+$ with nondecreasing c\`adl\`ag paths with $U_0=0$ and $U_{t}> 0$ for all $t>0$.
Let $\Omega$ be an open set in $\R^d$.
\begin{enumerate}[\ (a)]
\item The mapping $(x,u)\mapsto \P^W_x(\tau^W_\Omega>u)$ is $\mathcal{B}(\R^d)\times \mathcal{B}(\R_+)$-measurable. 
\item The mapping $x\mapsto \P_x(\tau^W_\Omega>U_t)$ is $\mathcal{B}(\R^d)$-measurable. 
\item The spectral heat content of the subordinate killed Brownian motion $(W^\Omega_{U_t})$ is well-defined. 
\item The regular heat content of the subordinate Brownian motion $(W_{U_t})$ is well-defined. 
\end{enumerate}
\end{lemma}

\pf
(a) Note that $\tau^W_\Omega>u$ if and only if $W_s\in\Omega$ for all $s\in[0,u]$, which implies $F_u:=\{\tau^W_\Omega>u\}\in \mathcal{B}(C[0,u])$. With this in mind, we prove the more general result that the mapping $(x,u)\mapsto \P^W_x(F_u)$ for a general $F_u\in \mathcal{B}(C[0,u])$ is $\mathcal{B}(\R^d)\times \mathcal{B}(\R_+)$-measurable. Since the family of sets $F_u\in \mathcal{B}(C[0,u])$ for which the mapping $(x,u)\mapsto \P^W_x(F_u)$ is $\mathcal{B}(\R^d)\times \mathcal{B}(\R_+)$-measurable forms a Dynkin system and the $\sigma$-algebra $\mathcal{B}(C[0,u])$ is generated by all finite-dimensional cylinder sets of the form $F_u=\{\omega\in C[0,u]: \omega(t_i)\in A_i \ \textrm{for} \ i=0,1,\ldots,k\}$, where $k\in\mathbb{N}$, $0=t_0< t_1<\cdots<t_k=u$, and  $A_i\in \mathcal{B}(\R)$ for $i=0,1,\ldots,k$, it suffices to prove the measurability for such cylinder sets only due to the Dynkin system theorem (\cite[Chapter 2, Theorem 1.3]{KaratzasShreve}). Now, for a cylinder set $F_u$ of the above form, 
\[
	\P^W_x(F_u)=\mathbf{1}_{A_0}(x)\int_{A_1}\int_{A_2}\cdots \int_{A_k} p(t_1;x,y_1)p(t_2-t_1;y_1,y_2)\cdots  p(t_k-t_{k-1};y_{k-1},y_k) dy_k\cdots dy_2 dy_1, 
\]
where $p(t;x,y):=(4\pi t)^{-d/2} e^{-\|x-y\|^2/(4t)}$ is the heat kernel. Since $t_k=u$, the mapping 
\[
	(x,u,y_1,\ldots,y_k)\mapsto \mathbf{1}_{A_0}(x)p(t_1;x,y_1)p(t_2-t_1;y_1,y_2)\cdots  p(t_k-t_{k-1};y_{k-1},y_k)
\]
 is $\mathcal{B}(\R^d)\times \mathcal{B}(\R_+)\times \mathcal{B}(\R^k)$-measurable. The desired measurability of $(x,u)\mapsto \P^W_x(F_u)$ is part of the consequence of the Fubini theorem. 

(b) Let $\FF^U$ denote the $\sigma$-algebra generated by $(U_t)$. Statement (a) together with the fact that the mapping $(x,\omega_2)\mapsto (x,U_t(\omega_2))$ is $\BB(\R^d)\times \FF^U$-measurable implies that the mapping $(x,\omega_2)\mapsto \P^W_x(\tau^W_\Omega>U_t(\omega_2))$ is $\BB(\R^d)\times \FF^U$-measurable. Thus, the Fubini theorem concludes that the mapping 
$
	x\mapsto \E^U[\P^W_x(\tau^W_\Omega>U_t)]= \P_x(\tau^W_\Omega>U_t) 
$
is $\BB(\R^d)$-measurable. 

(c) Note that $\{\tau_\Omega^{W^\Omega\circ U}> t\}=\{\tau_\Omega^W > U_t\}$ due to the assumption on the paths of $(U_t)$. This, together with statement (b), concludes that the mapping $x\mapsto \P_x(\tau_\Omega^{W^\Omega\circ U}> t)$ is $\mathcal{B}(\R^d)$-measurable. Hence, the integral $\int_\Omega \P_x(\tau_\Omega^{W^\Omega\circ U}> t) dx$ is well-defined, thereby yielding statement (c). 

(d) Since $\{W_u\in \Omega\}\in \BB(C[0,u])$, as in the arguments given in the proof of statements (a) and (b), one can show that the mapping $(x,u)\mapsto \P^W_x(W_u\in \Omega)$ is $\mathcal{B}(\R^d)\times \mathcal{B}(\R_+)$-measurable and consequently that the mapping $x\mapsto \E^U[\P^W_x(W_{U_t}\in \Omega)]=\P_x(W_{U_t}\in \Omega)$ is $\BB(\R^d)$-measurable. 
Hence, the integral $\int_\Omega \P_x(W_{U_t}\in \Omega) dx$ is well-defined, thereby yielding statement (d). 
\qed

\begin{prop}\label{prop:general_timechange}
Let $(W_t)$ be a Brownian motion in $\R^d$ independent of a process $(U_t)$ in $\R_+$ with nondecreasing c\`adl\`ag paths with $U_0=0$ and $U_{t}> 0$ for all $t>0$.
Let $\Omega$ be a bounded open interval when $d=1$, or a bounded connected $C^{1,1}$ open set when $d\geq 2$. 
Suppose that for any $\delta>0$,
\begin{align}\label{0009}
	\P^U(U_t>\delta)=o(\E^U[\sqrt{U_t}\mathbf{1}_{\{U_t\le \delta\}}]) \ \ \textrm{as} \ \ t\downarrow 0.
\end{align}
Suppose further that for any $\delta_{1}, \delta_{2}>0$, 
\begin{align}\label{eqn:independent limit}
\lim_{t\downarrow 0}\frac{\E^U[\sqrt{U_{t}}\mathbf{1}_{\{U_{t}\leq \delta_{1}\}}]}{\E^U[\sqrt{U_{t}}\mathbf{1}_{\{U_{t}\leq \delta_{2}\}}]}=1.
\end{align}
Then, for any $\delta>0$, 
\begin{align}\label{0010}
	\lim_{t\downarrow 0}\frac{|\Omega|-\tilde{Q}^{W\circ U}_{\Omega}(t)}{\E^U[\sqrt{U_{t}}\mathbf{1}_{\{U_{t}\le \delta\}}]}=\frac{2|\partial\Omega|}{\sqrt{\pi}}
\end{align}
and 
\begin{align}\label{0011}
	\lim_{t\downarrow 0}\frac{|\Omega|-\HH^{W\circ U}_{\Omega}(t)}{\E^U[\sqrt{U_t}\mathbf{1}_{\{U_t\le \delta\}}]}
	= \frac{|\partial\Omega|}{\sqrt{\pi}}.
\end{align}

\end{prop}
\pf
For any $\eps>0$, it follows from  \eqref{0031} that 
there exists $u_0>0$ such that
\[
	\left|\frac{|\Omega|-Q_{\Omega}^{W}(u)}{\sqrt{u}} -\frac{2|\partial \Omega|}{\sqrt{\pi}} \right| \le \eps 
	\text{ \ for all \ } 0<u\le u_0.
\]
Using the relation $\{\tau_\Omega^{W^\Omega\circ U}> t\}=\{\tau_\Omega^W > U_t\}$ and the Fubini theorem, for any $\delta>0$,
\begin{align*}
	&|\Omega|-\tilde{Q}^{W\circ U}_{\Omega}(t)
	=\int_\Omega\P^W_x\times \P^U(\tau_\Omega^{W^\Omega\circ U}\le t)dx
	=\E^U\left[\int_\Omega\P^W_x(\tau_\Omega^{W}\le U_t)dx\right]
	=\E^U[|\Omega|-Q^{W}_{\Omega}(U_t)]\\
	&=\E^U\Bigl[\frac{|\Omega|-Q^{W}_{\Omega}(U_t)}{\sqrt{U_t}}\sqrt{U_t}\mathbf{1}_{\{ U_t\le \delta\}}\Bigr]
		+\E^U[(|\Omega|-Q^{W}_{\Omega}(U_t))\mathbf{1}_{\{U_t> \delta\}}]\\
	&=:I_1(t)+I_2(t). 
\end{align*}
By assumption \eqref{0009}, 
\begin{align}\label{0036}
	0\le I_2(t)\le |\Omega|\P^U(U_t>\delta)=o(\E^U[\sqrt{U_t}\mathbf{1}_{\{U_t\le \delta\}}]) \ \ \textrm{as} \ \ t\downarrow 0. 
\end{align}
On the other hand, for any $0<\delta_{1}\leq u_{0}$,  
\[
	\left( \frac{2|\partial\Omega|}{\sqrt{\pi}}-\eps \right)\E^U[\sqrt{U_t}\mathbf{1}_{\{U_t\le \delta_1\}}] 
	\le I_1(t)
	\le \left( \frac{2|\partial\Omega|}{\sqrt{\pi}}+\eps \right)\E^U[\sqrt{U_t}\mathbf{1}_{\{U_t\le \delta_1\}}]. 
\]
This, together with \eqref{eqn:independent limit} and \eqref{0036}, implies that for any $\delta>0$, 
\begin{align*}
	 \frac{2|\partial\Omega|}{\sqrt{\pi}}-\eps
	 \le \liminf_{t\downarrow 0}\frac{|\Omega|-\tilde{Q}^{W\circ U}_{\Omega}(t)}{\E^U[\sqrt{U_t}\mathbf{1}_{\{U_t\le \delta\}}]}
	 \le \limsup_{t\downarrow 0}\frac{|\Omega|-\tilde{Q}^{W\circ U}_{\Omega}(t)}{\E^U[\sqrt{U_t}\mathbf{1}_{\{U_t\le \delta\}}]}
	\le \frac{2|\partial\Omega|}{\sqrt{\pi}}+\eps.
\end{align*}
Letting $\eps \downarrow 0$ yields \eqref{0010}, as desired.

Now, to establish \eqref{0011}, observe that
\begin{align*}
	&|\Omega|-\HH^{W\circ U}_{\Omega}(t)
	=|\Omega|-\int_\Omega\P^W_x\times \P^U(W_{U_t}\in\Omega)dx
	=|\Omega|-\E^U\left[\int_\Omega \P^W_x(W_{U_t}\in\Omega)dx\right]\\
	&=\E^U[|\Omega|-\HH^{W}_{\Omega}(U_t)]
	=\E^U\Bigl[\frac{|\Omega|-\HH^{W}_{\Omega}(U_t)}{\sqrt{U_t}}\sqrt{U_t}\mathbf{1}_{\{U_t\le \delta\}}\Bigr]
		+\E^U[(|\Omega|-\HH^{W}_{\Omega}(U_t))\mathbf{1}_{\{U_t> \delta\}}]. 
\end{align*}
Statement \eqref{0011} follows by \eqref{0032} together with an argument similar to 
that for the spectral heat content.
\qed

In the remainder of the paper, except in the proofs of Lemma \ref{lemma:continuous} and Theorem \ref{thm:complete expansion}, we use the simplified notations $\P$ and $\E$ in place of $\P^U$ and $\E^U$, respectively; this should not cause any confusion.

Sections \ref{section:SKBM} and \ref{section:inverse subordinators} investigate the regular and spectral heat contents for time-changed Brownian motions $(W_{U_t})$ with the independent time change $(U_t)$ being a subordinator $(D_t)_{t\ge 0}$ and its inverse $(E_t)_{t\ge 0}$, respectively, where 
$
	E_t:=\inf\{ u>0: D_u>t\}. 
$ 
Here, we briefly give a review of the concept of subordinators. 
By a subordinator $(D_t)$ starting at 0 with Laplace exponent $\phi$ with killing rate 0, drift 0, and infinite L\'evy measure $\nu$, we mean a one-dimensional strictly increasing L\'evy process with c\`adl\`ag paths starting at 0 with Laplace transform 
\[
	\E[e^{-sD_t}]=e^{-t\phi(s)}, \ \ \textrm{where} \  \phi(s)=\int_0^\infty (1-e^{-sy})\,\nu(dy), 
\]
with $\nu$ satisfying the conditions that $\nu((0,\infty))=\infty$ and $\int_0^\infty (y\wedge 1)\, \nu(dy)<\infty$. The distribution of $(D_t)$ is characterized by the Laplace exponent $\phi$, which is a Bernstein function on $(0,\infty)$ 
(i.e.,\ $\phi\in C^\infty(0,\infty)$ with $\phi\ge 0$ and $(-1)^n\phi^{(n)}\le 0$ for all $n\in\N$) with $\phi(0^+)=0$. 
In particular, since $\phi(s_2)-\phi(s_1)=\int_{0}^{\infty}(e^{-s_1 y}-e^{-s_2 y})\nu(dy)>0$ for $0<s_1<s_2$, $\phi$ is strictly increasing and concave. 
Moreover, the assumption $\nu((0,\infty))=\infty$ implies $\phi(s)\to \infty$ as $s\to\infty$.
For a general account of subordinators, see \cite{Bertoin_Levy}.

To describe the wide class of Laplace exponents to be considered in this paper, we recall the notion of regular variation. By definition, a function $f:(0,\infty)\to (0,\infty)$ is regularly varying at $\infty$ (resp.\ at $0$) with index $\rho\in\R$ if $f(cs)/f(s)\to c^\rho$ as $s\to \infty$ (resp.\ as $s\to 0$) for each fixed $c>0$, or equivalently, $f$ has representation $f(s)=s^\rho \ell(s)$ with some slowly varying function $\ell$ at $\infty$ (resp.\ at $0$); 
i.e.,\ $\ell(cs)/\ell(s)\to 1$ as $s\to\infty$ (resp.\ as $s\to 0$) for each fixed $c>0$. 
We use the notation $\mathcal{R}_\rho(\infty)$ (resp.\ $\mathcal{R}_\rho(0^+)$) to denote the class of all regularly varying functions at $\infty$ (resp.\ at 0) with index $\rho\in \R$. Note that for any $f \in \mathcal{R}_\rho(\infty)$, as $s\to \infty$, $f(s)\to \infty$ if $\rho>0$ and $f(s)\to 0$ if $\rho<0$ (see \cite[Proposition 1.5.1]{Bingham_book}). 

This paper deals with a subordinator $(D_t)$ with Laplace exponent $\phi\in \mathcal{R}_\beta(\infty)$, where $\beta\in(0,1)$. Examples of such Laplace exponents include $\phi(s)=s^\beta$ for a $\beta$-stable subordinator and $\phi(s)=(s+\theta)^\beta-\theta^\beta$ for a tempered stable subordinator with tempering factor $\theta>0$. 
From theoretical viewpoints, \textit{with the use of the general form of the Laplace exponent, the scaling property (i.e.,\ self-similarity) that a stable subordinator $(S_{t}^{(\beta)})$ possesses (and hence so does its inverse $(E_t^{(\beta)})$) is no longer available;} therefore, one often needs to carefully examine asymptotic behaviors of various quantities using different theoretical tools and ideas.

\section{Time change by subordinators}\label{section:SKBM}
In this section, we apply Proposition \ref{prop:general_timechange} to study the asymptotic limit of the spectral heat content for time-changed Brownian motions when the underlying time change $(U_t)$ is given by a subordinator $(D_t)$ with Laplace exponent $\phi\in \mathcal{R}_\beta(\infty)$, where $\beta\in(0,1)$.
This application recovers and generalizes a result on symmetric stable processes established in \cite[Theorem 1.1]{ParkSong19}. 
The section is divided into the following three separate cases: $\beta\in(\frac 12,1)$, $\beta=\frac 12$, and $\beta\in(0,\frac 12)$. 
The main results corresponding to the three cases appear in Theorems \ref{thm:stable2}, \ref{thm:withlowerorder}, and \ref{thm:lessthan1/2}, respectively. 

We first establish an upper bound for the heat kernel for $(D_{t})$, which is global in space and is valid for all $\phi\in \mathcal{R}_{\beta}(\infty)$ with $\beta\in (0,1)$.
The result essentially follows from \cite{Grzywny2018}, but we state it here in a form that is convenient for our purpose. 

\begin{lemma}\label{lemma:HK ub}
Suppose $(D_t)$ is a subordinator with Laplace exponent $\phi\in \mathcal{R}_{\beta}(\infty)$ with $\beta\in (0,1)$ and its L\'evy measure has an almost monotone density.
Then, the transition density $p(t,x)$ of $(D_{t})$ exists. 
Moreover, there exist $t_{0}>0$ and $c>0$ such that 
\beq\label{eqn:HK ub}
p(t,x)\leq c \left(\phi^{-1}(1/t)\wedge tx^{-1}\phi(1/x)\right)
\eeq
for all $0<t<t_0$ and $x>0$. 
\end{lemma}

\begin{remark}
\begin{em}
(a) Lemma \ref{lemma:HK ub} is deeply connected with \cite[Theorem 4.8]{Grzywny2018}, which provides an estimate similar to \eqref{eqn:HK ub} for all $x$ away from 0.
The proof below clarifies how to derive \eqref{eqn:HK ub} for all $x>0$ (including values near 0) based on some results in \cite{Grzywny2018}.

(b) A number of subordinators have Laplace exponents that are \textit{complete} Bernstein functions (see e.g.\ \cite[Chapter 16]{Schilling}), and hence, each of their L\'evy measures has an almost monotone density.

\end{em}
\end{remark}

\pf
Suppose $\phi\in \mathcal{R}_{\beta}(\infty)$ with $\beta\in(0,1)$, and fix a constant $\eta$ such that $0<\beta-\eta<\beta+\eta<1$. 
By Potter's Theorem \cite[Theorem 1.5.6. (iii)]{Bingham_book}, there exists a constant $x_0>0$ such that 
$$
\frac{\phi(y)}{\phi(z)}\leq 2\max\left(\left(\frac{y}{z}\right)^{\beta-\eta}, \left(\frac{y}{z}\right)^{\beta+\eta}\right)  \text{ for all } y,z\geq x_0.
$$
Writing $z=\lam y$ with $\lam\ge 1$ yields the weak lower scaling condition at $\infty$ given by
\beq\label{eqn:scaling_WLSC}
\frac{\phi(\lam y)}{\phi(y)}\geq \frac12 \lam^{\beta-\eta} \text{ for all } y\geq x_0, \lam\geq 1.
\eeq
On the other hand, setting $y=\lam z$ with $\lam\ge 1$ gives the weak upper scaling condition 
\beq\label{eqn:scaling_WUSC}
\frac{\phi(\lam z)}{\phi(z)}\leq 2 \lam^{\beta+\eta} \text{ for all } z\geq x_0, \lam\geq 1. 
\eeq
In the notations used in \cite{Grzywny2018}, conditions \eqref{eqn:scaling_WLSC} and \eqref{eqn:scaling_WUSC} are expressed as $\phi\in \text{WLSC}(\beta -\eta,\frac12, x_{0})\cap \text{WUSC}(\beta+\eta, 2, x_{0})$. This implies $-\phi^{''}\in \text{WLSC}(\beta -\eta-2,C_1, x_{0})\cap \text{WUSC}(\beta+\eta-2, C_2, x_{0})$ for some $0<C_1< 1< C_2$, and consequently, the transition density $p(t,x)$ of $(D_t)$ exists (see \cite[Corollary 2.7, Theorem 3.3]{Grzywny2018} for details).

It follows from \cite[Equations (4.13) and (4.14)]{Grzywny2018} and \cite[Equation (4.29)]{Grzywny2018} that there exists $t_{0}>0$ such that
$$
p(t,x)\leq c\left(\phi^{-1}(1/t)\wedge tx^{-1}\phi(1/x)\right)
$$
for all $0<t\leq t_0$ and $x\geq 2et\phi'(\psi^{-1}(1/t))$, 
where $\psi$ is the characteristic exponent of $(D_t)$ 
and 
$\psi^{-1}$ is the inverse of the symmetric, continuous, and nondecreasing majorant of $\psi^{*}(r)=\sup_{|z|\leq r}\mathfrak{R}\psi(z)$, with $\mathfrak{R}(\cdot)$ denoting the real part of the argument.
Note that for all sufficiently small $t$, it follows from  \cite[Proposition 2.3]{Grzywny2018} and \cite[Equation (4.10)]{Grzywny2018} that
$$
t\phi'(\psi^{-1}(1/t))\approx t\psi^{-1}(1/t)^{-1}\phi(\psi^{-1}(1/t))\approx \phi^{-1}(1/t)^{-1};
$$
hence, for $x< 2et\phi'(\psi^{-1}(1/t))$ we have $\phi(1/x)\gtrsim1/t$ and $\frac{t}{x}\phi(1/x)\gtrsim \frac{1}{x} \gtrsim \phi^{-1}(1/t)$, where 
$A(x)\gtrsim B(x)$ means there exists a constant $c$ such that $A(x) \geq c B(x)$ and $A(x)\approx B(x)$ means $A(x)\gtrsim B(x)$ and $B(x)\gtrsim A(x)$ on their domains.
Hence, if $x< 2et\phi'(\psi^{-1}(1/t))$, then $\frac{t}{x}\phi(1/x)\gtrsim \phi^{-1}(1/t)$, and it follows from \cite[Lemma 6]{KS15} together with \cite[Equation (4.21)]{Grzywny2018} and \cite[Equations (4.2) and (4.13)]{Grzywny2018} that $p(t,x)\leq c\phi^{-1}(1/t)$.
This establishes \eqref{eqn:HK ub} for all small $t$ and $x>0$.
\qed

\subsection{Case 1: $\phi\in \mathcal{R}_\beta(\infty)$ with $\beta\in(\frac 12,1)$}
We first focus on the case when $\phi\in \mathcal{R}_\beta(\infty)$ with $\beta\in(\frac 12,1)$. 
In order to apply Proposition \ref{prop:general_timechange}, we must verify that conditions \eqref{0009} and \eqref{eqn:independent limit} hold. The first lemma concerns condition \eqref{0009}.
\begin{lemma}\label{lemma:stable1}
Suppose that $(D_{t})$ is a subordinator with Laplace exponent $\phi\in \mathcal{R}_{\beta}(\infty)$ with $\beta\in (\frac12,1)$. 
Then condition \eqref{0009} holds.
\end{lemma}

\pf
To verify condition \eqref{0009}, note first by \cite[Corollary 4.14]{Rosinski_isomorphism} that 
\begin{align}\label{0016}
	\P(D_t\ge \delta)\sim \nu([\delta,\infty)) t
	\ \ \textrm{as} \ \ t\downarrow 0.
\end{align}
(Also see \cite[Proposition 1]{Kobayashi_smallball}, which expresses \eqref{0016} in terms of the inverse $(E_t)$ of $(D_t)$ under the assumption that $\nu(\{\delta\})=0$.) 
Next, recalling that $\phi$ is continuous and strictly increasing and hence has inverse $\phi^{-1}$, define a scaled process $(Y_s^{(t)})_{s\ge 0}$ for each fixed $t>0$ by 
\beq\label{eqn:approximate scaling}
Y_{s}^{(t)}:=\phi^{-1}(1/t)D_{ts}, \quad s\ge 0.
\eeq
Then the Laplace exponent $\phi^{(t)}(\lam)$ of $(Y_{s}^{(t)})_{s\ge 0}$ is given by $t\phi(\lam \phi^{-1}(1/t))$ and 
$$
\lim_{t\downarrow 0}\phi^{(t)}(\lam)=\lam^{\beta},
$$
as observed in the proof of \cite[Lemma 4.5]{GPS19}. (Note that the assumption of regular variation of $\phi$ is crucial here and cannot be replaced by the weak scaling conditions \eqref{eqn:scaling_WLSC}-\eqref{eqn:scaling_WUSC}.)
Hence, by \cite[Theorem XIII.1 2]{FellerTwo71}, 
\begin{align}\label{eqn:weak convergence}
	Y_{s}^{(t)}\to S_{s}^{(\beta)} \ \ \textrm{in distribution as} \ \ t\downarrow 0
\end{align}
 for each fixed $s\ge 0$,  
where $(S_{s}^{(\beta)})_{s\ge 0}$ is a $\beta$-stable subordinator. (The convergence actually holds for all finite-dimensional distributions due to the independent and stationary increments.) 

The condition that $\phi\in \mathcal{R}_\beta(\infty)$ implies that $\phi(\lambda)\to \infty$ as $\lambda\to \infty$, and 
its inverse satisfies $\phi^{-1} \in \mathcal{R}_{1/\beta}(\infty)$, which in turn implies $\phi^{-1}(1/t) \in \mathcal{R}_{-1/\beta}(0^+)$. 
In particular, $\phi^{-1}(1/t)\to \infty$ as $t\downarrow 0$, so for small enough $t>0$, the inequality $\delta \phi^{-1}(1/t)\ge 1$ holds, and hence,
\begin{align*}
\E[D_{t}^{1/2} \mathbf{1}_{\{D_{t}\leq \delta\}}]
&=\E[[\phi^{-1}(1/t)]^{-1/2}(Y_{1}^{(t)})^{1/2} \mathbf{1}_{\{[\phi^{-1}(1/t)]^{-1}Y_{1}^{(t)}\leq \delta\}}]\\
&\geq[\phi^{-1}(1/t)]^{-1/2}\E[(Y_{1}^{(t)})^{1/2} \mathbf{1}_{\{Y_{1}^{(t)}\leq 1\}}].
\end{align*}
This, together with \eqref{0016}, yields 
\begin{align*}
\limsup_{t\downarrow 0}\frac{\P(D_{t}>\delta)}{\E[\sqrt{D_{t}} 1_{\{D_{t}\leq \delta\}}]}
&\leq 2\nu([\delta,\infty))\limsup_{t\downarrow  0} \frac{t[\phi^{-1}(1/t)]^{1/2}}{\E\left[(Y_{1}^{(t)})^{1/2} \mathbf{1}_{\{Y_{1}^{(t)}\leq 1\}}\right]}. 
\end{align*}
The latter limit is 0 due to the weak convergence in \eqref{eqn:weak convergence} as well as the convergence $t[\phi^{-1}(1/t)]^{1/2}\to 0$, which follows from the fact that $\phi^{-1}(1/t) \in \mathcal{R}_{-1/\beta}(0^+)$ with $\beta\in(\frac 12,1)$. Thus, condition \eqref{0009} holds.
\qed

Proposition \ref{prop:general_timechange} requires the exact asymptotic rate of $\E[\sqrt{D_{t}} \mathbf{1}_{\{D_{t}\leq \delta\}}]$ as $t\downarrow 0$, which will be derived in Proposition \ref{lemma:convergence to SS}. The following technical lemma is needed to prove the proposition.
\begin{lemma}\label{lemma:technical1}
Suppose that $\phi\in \mathcal{R}_{\beta}(\infty)$ with $\beta\in (\frac12,1)$. 
Then for any $\eps>0$ and $\delta>0$, there exist $M=M(\eps,\beta)\ge 1$ and $t_{0}>0$ such that 
\beq\label{eqn:cond3}
M[\phi^{-1}(1/t)]^{-1}<\delta \ \ \textrm{and} \ \ t[\phi^{-1}(1/t)]^{1/2}\int_{M[\phi^{-1}(1/t)]^{-1}}^{\delta}u^{-1/2}\phi(1/u)du \leq \eps
\eeq
for all $0<t< t_{0}$.
\end{lemma}

\pf
Since $\phi\in\mathcal{R}_{\beta}(\infty)$ with $\beta\in (\frac12,1)$, we can take small $\eta>0$ so that $\gamma:=\beta-\eta>\frac 12$. Then the weak lower scaling condition \eqref{eqn:scaling_WLSC} reads
\beq\label{eqn:scaling}
\frac{\phi(\lam y)}{\phi(y)}\geq \frac12 \lam^{\gamma} \text{ for all } y\geq x_0, \lam\geq 1,
\eeq
where $x_0>0$ is some constant.

Let $\eps>0$ and $\delta>0$. Let $c_1:=\frac{2}{\beta-1/2}$. Take $M=M(\eps,\beta)\ge 1$ so that $2c_{1}M^{-(\gamma-\frac12)}<\eps$. 
Note that $\phi^{-1}(1/t)\to \infty$ as $t\downarrow 0$. Moreover, by Karamata's Tauberian Theorem (\cite[Theorem 1.5.11]{Bingham_book}), 
\[
	\int_{1/\delta}^{\phi^{-1}(1/t)/M}v^{-3/2}\phi(v)dv \sim \frac{1}{\beta-1/2}\left[\frac{\phi^{-1}(1/t)}{M}\right]^{-1/2}\phi\left(\frac{\phi^{-1}(1/t)}{M}\right),
\]
 so 
there exists $t_0=t_0(\eps,\beta,\delta)>0$ such that for all $0<t<t_0$, 
\begin{align}\label{0038}
	\frac{\phi^{-1}(1/t)}{M}>\max\left(x_0,\frac{1}{\delta}\right)
\end{align}
 and 
\begin{align}\label{0033}
\int_{M[\phi^{-1}(1/t)]^{-1}}^{\delta}u^{-1/2}\phi(1/u)du
=\int_{1/\delta}^{\phi^{-1}(1/t)/M}v^{-3/2}\phi(v)dv
\le c_{1}\left[\frac{\phi^{-1}(1/t)}{M}\right]^{-1/2}\phi\left(\frac{\phi^{-1}(1/t)}{M}\right).
\end{align}
For $0<t<t_0$, the first inequality in \eqref{eqn:cond3} follows immediately from \eqref{0038}. Moreover, by \eqref{eqn:scaling} and \eqref{0038},  
$$
\frac{1/t}{\phi\left(\frac{\phi^{-1}(1/t)}{M}\right)}=\frac{\phi\left(M\cdot \frac{\phi^{-1}(1/t)}{M}\right)}{\phi\left(\frac{\phi^{-1}(1/t)}{M}\right)}\geq \frac12 M^{\gamma}. 
$$
This together with \eqref{0033} gives
\begin{align*}
t[\phi^{-1}(1/t)]^{1/2}\int_{M[\phi^{-1}(1/t)]^{-1}}^{\delta}u^{-1/2}\phi(1/u)du
\leq c_{1}t \phi\left(\frac{\phi^{-1}(1/t)}{M}\right)M^{1/2}
\leq 2c_{1}M^{-(\gamma-\frac12)}<\eps, 
\end{align*}
yielding the second inequality in \eqref{eqn:cond3}, as desired.
\qed

It follows from \cite[Equation (2.5)]{Valverde2014} that a $\beta$-stable subordinator $(S_{t}^{(\beta)})_{t\ge 0}$ has moments of orders $\gamma<\beta$ given by the following explicit formula: 
\beq\label{eqn:SS moment}
\E[(S_{1}^{(\beta)})^{\gamma}]=\frac{\Gamma(1-\frac{\gamma}{\beta})}{\Gamma(1-\gamma)}, \quad -\infty <\gamma<\beta.
\eeq
In the following arguments, we often use the notation $\E[X,A]$ in place of $\E[X\mathbf{1}_A]$.

\begin{prop}\label{lemma:convergence to SS}
Suppose $(D_{t})$ is a subordinator with Laplace exponent $\phi\in \mathcal{R}_{\beta}(\infty)$ with $\beta \in (\frac12,1)$ and its L\'evy measure has an almost monotone density.
Then for any $\delta>0$, 
$$
\E[\sqrt{D_{t}}\mathbf{1}_{\{D_{t}\leq \delta\}}]\sim [\phi^{-1}(1/t)]^{-1/2}\E[(S_{1}^{(\beta)})^{1/2}] \text{ as } t\downarrow 0,
$$
where $(S_{t}^{(\beta)})_{t\ge 0}$ is a $\beta$-stable subordinator. 
In particular, condition \eqref{eqn:independent limit} holds.
\end{prop}
\pf
Define the scaled process $(Y_{s}^{(t)})_{s\ge 0}$ as in \eqref{eqn:approximate scaling}, which approximates in distribution the $\beta$-stable subordinator $(S_{s}^{(\beta)})_{s\ge 0}$ as stated in \eqref{eqn:weak convergence}. 
Let $\eps>0$. 
Since $\beta\in (\frac12,1)$, it follows from \eqref{eqn:SS moment} and the dominated convergence theorem that there exists $M_{1}>0$ such that
\begin{align}\label{0034}
\E[(S_{1}^{(\beta)})^{1/2}, S_{1}^{(\beta)}> M_{1}] <\eps.
\end{align}
By the proof of Lemma \ref{lemma:technical1}, we can take $M\ge M_1$ large enough and $t_0>0$ small enough so that condition \eqref{eqn:cond3} holds for all $0<t<t_0$.

Fix $0<t<t_0$. Noting that $M[\phi^{-1}(1/t)]^{-1}< \delta$, decompose $\E[\sqrt{D_{t}}\mathbf{1}_{\{D_{t}\leq \delta\}}]=\E[D_{t}^{1/2}, D_{t}\leq \delta]$ as
\beq\label{eqn:cases1}
\E[D_{t}^{1/2}, D_{t}\leq \delta]=\E[D_{t}^{1/2}, D_{t}\leq M[\phi^{-1}(1/t)]^{-1}]+ \E[D_{t}^{1/2}, M[\phi^{-1}(1/t)]^{-1} <D_{t}\leq \delta]. 
\eeq
Express the first term on the right hand side as
\begin{eqnarray}\label{eqn:cases2}
\E[D_{t}^{1/2}, D_{t}\leq M[\phi^{-1}(1/t)]^{-1}]
=[\phi^{-1}(1/t)]^{-1/2}\E[(Y_{1}^{(t)})^{1/2},  Y_{1}^{(t)} \leq M].
\end{eqnarray}
Dropping the second term on the right hand side of \eqref{eqn:cases1}, combining \eqref{eqn:cases2} with \eqref{eqn:weak convergence}, and using \eqref{0034}, we obtain
\begin{align*}
\liminf_{t\downarrow 0}[\phi^{-1}(1/t)]^{1/2}\E[D_{t}^{1/2}, D_{t}\leq \delta] 
&\geq\liminf_{t\downarrow 0}[\phi^{-1}(1/t)]^{1/2}\E[D_{t}^{1/2}, D_{t}\leq M[\phi^{-1}(1/t)]^{-1}] \\
&\geq \E[(S_{1}^{(\beta)})^{1/2}, S_{1}^{(\beta)}\leq M] \geq \E[(S_{1}^{(\beta)})^{1/2}] -\eps.
\end{align*}
Since $\eps>0$ is arbitrary, the latter yields the lower bound  
\beq\label{eqn:lb1}
\liminf_{t\downarrow 0}[\phi^{-1}(1/t)]^{1/2}\E[D_{t}^{1/2}, D_{t}\leq \delta] \geq \E[(S_{1}^{(\beta)})^{1/2}].
\eeq

Now we derive the upper bound. 
By Lemma \ref{lemma:HK ub}, there exist $t_{0}>0$ and $c>0$ such that 
$
	p(t,u)\leq ctu^{-1}\phi(1/u)
$
for all $0<t<t_0$ and $u>0$.
 Hence, the second term on the right hand side of \eqref{eqn:cases1}, which equals $\int_{M[\phi^{-1}(1/t)]^{-1}}^{\delta}u^{1/2}p(t,u)du$,  is bounded above by
\begin{eqnarray*}
\int_{M[\phi^{-1}(1/t)]^{-1}}^{\delta}cu^{1/2}\frac{t}{u}\phi(1/u)du=ct\int_{M[\phi^{-1}(1/t)]^{-1}}^{\delta}u^{-1/2}\phi(1/u)du.
\end{eqnarray*}
This, together with \eqref{eqn:cond3}, yields
\begin{align}\label{eqn:ub1}
[\phi^{-1}(1/t)]^{1/2}\E[D_{t}^{1/2}, M[\phi^{-1}(1/t)]^{-1} <D_{t}\leq \delta]
\leq ct[\phi^{-1}(1/t)]^{1/2}\int_{M[\phi^{-1}(1/t)]^{-1}}^{\delta}u^{-1/2}\phi(1/u)du
\leq c\eps
\end{align}
for all $0<t\leq t_{0}$.
Combining \eqref{eqn:weak convergence}, \eqref{eqn:cases1}, \eqref{eqn:cases2}, and \eqref{eqn:ub1} gives
\begin{align*}
\limsup_{t\downarrow 0}[\phi^{-1}(1/t)]^{1/2}\E[D_{t}^{1/2}, D_{t}\leq \delta] 
\leq \E[(S_{1}^{(\beta)})^{1/2}, (S_{1}^{(\beta)})^{1/2}\leq M ]+c\eps
\leq \E[(S_{1}^{(\beta)})^{1/2} ]+c\eps.
\end{align*}
Since $\eps>0$ is arbitrary, the latter yields the upper bound
\beq\label{eqn:ub2}
\limsup_{t\downarrow 0}[\phi^{-1}(1/t)]^{1/2}\E[D_{t}^{1/2}, D_{t}\leq \delta] \leq \E[(S_{1}^{(\beta)})^{1/2}].
\eeq
The desired conclusion now follows from \eqref{eqn:lb1} and \eqref{eqn:ub2}.
\qed

Note that if $(D_t)$ itself is a $\beta$-stable subordinator $(S^{(\beta)}_t)$, then the asymptotic relation in Proposition \ref{lemma:convergence to SS} follows immediately from the self-similarity of $(D_t)$ with index $1/\beta$. 
Combining Lemma \ref{lemma:stable1} and Proposition \ref{lemma:convergence to SS} gives the following theorem for the case when $\beta\in (\frac12,1)$. 
\begin{thm}\label{thm:stable2}
Let $\Omega$ be a bounded open interval when $d=1$, or a bounded connected $C^{1,1}$ open set when $d\geq 2$.
Suppose $(D_{t})$ is a subordinator with Laplace exponent $\phi\in \mathcal{R}_{\beta}(\infty)$ with $\beta \in (\frac12,1)$ and its L\'evy measure has an almost monotone density. Let $(W_t)$ be a Brownian motion independent of $(D_t)$.
Then
\begin{align}\label{0050}
\lim_{t\downarrow 0} \frac{|\Omega|-\tilde{Q}_{\Omega}^{W\circ D}(t)}{[\phi^{-1}(1/t)]^{-1/2}}=\E[(S_{1}^{(\beta)})^{1/2}]\times \frac{2|\partial \Omega|}{\sqrt{\pi}},
\end{align}
and 
\begin{align}\label{0051}
\lim_{t\downarrow 0} \frac{|\Omega|-\HH_{\Omega}^{W\circ D}(t)}{[\phi^{-1}(1/t)]^{-1/2}}=\E[(S_{1}^{(\beta)})^{1/2}]\times \frac{|\partial \Omega|}{\sqrt{\pi}},
\end{align}
where $(S_{t}^{(\beta)})_{t\ge 0}$ is a $\beta$-stable subordinator. 
\end{thm}

\begin{remark}\label{remark:interpretation1}
\begin{em}
It is possible to rewrite \eqref{0050} and \eqref{0051} respectively as follows: 
\begin{align}\label{0052}
	\lim_{t\downarrow 0} \frac{|\Omega|-\tilde{Q}_{\Omega}^{W\circ D}(t)}{[\phi^{-1}(1/t)]^{-1/2}}=\E[\sup_{0\le u\leq S_{1}^{(\beta)}}B_{u}]|\partial \Omega|
 \ \ \textrm{and} \ \ 
\lim_{t\downarrow 0} \frac{|\Omega|-\HH_{\Omega}^{W\circ D}(t)}{[\phi^{-1}(1/t)]^{-1/2}}=\frac12\E[\sup_{0\le u\leq S_{1}^{(\beta)}}B_{u}]|\partial \Omega|,
\end{align}
where $(B_t)$ is a \textit{one-dimensional} Brownian motion starting at 0 with $\E[e^{i\xi B_t}]=e^{-t\xi^2}$, and $(S_{t}^{(\beta)})$ is a $\beta$-stable subordinator independent of $(B_t)$. Indeed, by \cite[Problem 2.8.2]{KaratzasShreve}, the first moment of the running maximum of the Brownian motion is given by   
\begin{align}\label{0054}
	\E[\sup_{0\le u\le t}B_u]
	=\int_0^\infty \frac{x}{\sqrt{\pi t}}e^{-x^2/(4t)}dx
	=\frac{2\sqrt{t}}{\sqrt{\pi}}, \ \ t>0, 
\end{align}
which, together with the independence assumption, yields $\E[\sup_{0\le u\leq S_{1}^{(\beta)}}B_{u}]=\E[(S_{1}^{(\beta)})^{1/2}]\times \frac{2}{\sqrt{\pi}}$. 
The asymptotic limit of the spectral heat content given in \eqref{0052} can be interpreted as follows: at every point near the boundary, only the one dimensional shortest path, which is parallel to the orthogonal vector at the boundary point, contributes to the spectral heat content.
\end{em}
\end{remark}

\subsection{Case 2: $\phi(\lam)=\lam^{1/2}+\phi_{1}(\lam)$ with $\phi_{1}\in \mathcal{R}_{\beta_1}(\infty)$ with $\beta_1\in (0,1/2)$}
We now turn our attention to the case when $\phi\in \mathcal{R}_{\beta}(\infty)$ with $\beta=1/2$.  
The simplest case is when $\phi(\lambda)=\lambda^{1/2}$, with which $(D_t)$ is the L\'evy subordinator and the time-changed Brownian motion $(W_{D_t})$ has the Cauchy distribution. 
The following result was already established in \cite[Theorem 1.1]{ParkSong19} and is now recovered via Proposition \ref{prop:general_timechange}. 
\begin{prop}\label{prop:stable2}
Let $(W_t)$ be a Brownian motion independent of a L\'evy subordinator $(S^{(1/2)}_t)$ whose Laplace exponent is given by $\phi(\lambda)=\lambda^{1/2}$. Let $\Omega$ be a bounded open interval when $d=1$, or a bounded connected $C^{1,1}$ open set when $d\geq 2$. Then  
\begin{align}\label{0015}
	\lim_{t\downarrow 0}\frac{|\Omega|-\tilde{Q}^{W\circ S^{(1/2)}}_{\Omega}(t)}{t\log(1/t)}=\frac{2|\partial\Omega|}{\pi}
\ \  \ \textrm{and} \ \ \ 
	\lim_{t\downarrow 0}\frac{|\Omega|-\HH^{W\circ S^{(1/2)}_{\Omega}(t)}}{t\log(1/t)}= \frac{|\partial\Omega|}{\pi}. 
\end{align}
\end{prop}
\pf
For a fixed $\delta>0$, by the self-similarity of $(S^{(1/2)}_t)$ with index $\frac 1\beta=2$ as well as \cite[Lemma 3.2]{ParkSong19}, 
\[
	\E[\sqrt{S^{(1/2)}_t}\mathbf{1}_{\{S^{(1/2)}_t\le \delta\}}]
	=t\E[\sqrt{S^{(1/2)}_1}\mathbf{1}_{\{S^{(1/2)}_1\le \delta t^{-2}\}}]
	\sim \frac{t\log (1/t)}{\sqrt{\pi}} \ \ \textrm{as} \ \ t\downarrow 0.
\]
Hence, condition \eqref{eqn:independent limit} holds.
This, together with \eqref{0016}, gives
\[
	\frac{\P(S^{(1/2)}_t>\delta)}{\E[\sqrt{S^{(1/2)}_t}\mathbf{1}_{\{S^{(1/2)}_t\le \delta\}}]}
	\sim \frac{\sqrt{\pi}\nu([\delta,\infty))}{\log(1/t)}
	\ \ \textrm{as} \ \ t\downarrow 0, 
\]
so condition \eqref{0009} holds. Application of Proposition \ref{prop:general_timechange} now yields \eqref{0015}, as desired. 
\qed

Proposition \ref{prop:stable2} can be described in a more general setting where 
$\phi(\lam)=\lam^{1/2} + $ any lower order term.
We start with the following lemma, which is an analogue of \cite[Lemma 3.2]{ParkSong19}. 

\begin{lemma}\label{lemma:mixedconv}
Let $\delta>0$ be a fixed constant. 
Suppose that the Laplace exponent $\phi$ of a subordinator $(D_{t})$ takes the form $\phi(\lam)=\lam^{1/2} +\phi_{1}(\lam)$, where $\phi_1$ is 
the Laplace exponent of a subordinator whose L\'evy measure has an almost monotone density and
$\phi_{1} \in\mathcal{R}_{\beta_1}(\infty)$ with $\beta_1\in (0,\frac12)$.
Then
$$
\lim_{t\downarrow 0}\frac{\E[\sqrt{D_t}\mathbf{1}_{\{D_t\le \delta\}}]}{t\log(1/t)}=\frac{1}{\sqrt{\pi}}
$$
\end{lemma}

\pf
Express $(D_{t})$ as the sum of independent subordinators
$$
D_{t}=S_{t}^{(1/2)} +D_{t}^{\phi_1},
$$
where $(S_{t}^{(1/2)})$ is L\'evy subordinator and $(D_{t}^{\phi_1})$ is a subordinator with the Laplace exponent $\phi_1$. 
Note that 
\begin{align}\label{eqn:aux-1}
\E[\sqrt{D_t}\mathbf{1}_{\{D_t\le \delta\}}]=\E[\left(S_{t}^{(1/2)} + D_{t}^{\phi_1}\right)^{1/2}, S_{t}^{(1/2)}+D_{t}^{\phi_1}\le \delta].
\end{align}
Since the function $x\mapsto \sqrt{x}$ is subadditive, \eqref{eqn:aux-1} is bounded above by
\begin{align*}
&\E[\left(S_{t}^{(1/2)}\right)^{1/2}, S_{t}^{(1/2)}+D_{t}^{\phi_{1}}\le \delta]
+\E[\left( D_{t}^{\phi_{1}}\right)^{1/2}, S_{t}^{(1/2)}+D_{t}^{\phi_{1}}\le \delta]\\
&\leq\E[\left(S_{t}^{(1/2)}\right)^{1/2}, S_{t}^{(1/2)}\le \delta]+\E[\left(D_{t}^{\phi_{1}}\right)^{1/2}, D_{t}^{\phi_{1}}\le \delta].
\end{align*}
By Lemma \ref{lemma:HK ub}, there exist $t_{0}>0$ and $c_1>0$ such that  the transition density $p^{\phi_{1}}(t,u)$ of $(D_{t}^{\phi_{1}})$ satisfies $p^{\phi_1}(t,u)\leq c_1 tu^{-1}\phi_{1}(1/u)$ for all $0<t<t_0$ and $u>0$.
Hence, for $0<t<t_0$, 
\begin{align*}
\E[\left(D_{t}^{\phi_{1}}\right)^{1/2}, D_{t}^{\phi_{1}}\le \delta]
\leq c_1\int_{0}^{\delta}u^{1/2}\frac{t}{u}\phi_{1}(1/u)du
=c_{1}t\int_{1/\delta}^{\infty}v^{-3/2}\phi_{1}(v)dv.
\end{align*}
Note that the integral $\int_{1/\delta}^{\infty}v^{-3/2}\phi_{1}(v)dv$ is finite since $\phi_{1}\in \mathcal{R}_{\beta_1}(\infty)$ with $\beta_1\in (0,\frac12)$ (see \cite[Proposition 1.5.10]{Bingham_book}). Hence, 
$$
\lim_{t\downarrow 0}\frac{\E[\left(D_{t}^{\phi_{1}}\right)^{1/2}, D_{t}^{\phi_{1}}\le \delta]}{t\log(1/t)}=0.
$$
On the other hand, by \cite[Lemma 3.2]{ParkSong19}, 
$$
\lim_{t\downarrow 0}\frac{\E[\left(S_{t}^{(1/2)}\right)^{1/2}, S_{t}^{(1/2)}\le \delta]}{t\log(1/t)}=\frac{1}{\sqrt{\pi}}.
$$
Hence, we have shown the upper bound 
$$
\limsup_{t\downarrow 0}\frac{\E[\sqrt{D_t}\mathbf{1}_{ \{ D_t\le \delta\}  } ]}{t\log(1/t)}\leq \frac{1}{\sqrt{\pi}}.
$$
In terms of the lower bound, by the independence between $(S^{1/2}_t)$ and $(D^{\phi_1}_t)$, 
\begin{align*}
&\E[\left(S_{t}^{(1/2)} + D_{t}^{\phi_{1}}\right)^{1/2}, S_{t}^{(1/2)}+D_{t}^{\phi_{1}}\le \delta]
\geq \E[\left( S_{t}^{(1/2)}\right)^{1/2}, S_{t}^{(1/2)}+D_{t}^{\phi_{1}}\le \delta]\\
=&\int_{0}^{\delta}\E[\left( S_{t}^{(1/2)}\right)^{1/2}, S_{t}^{(1/2)}\le \delta-y]\P(D_{t}^{\phi_{1}}\in dy)\geq \int_{0}^{\delta/2}\E[\left( S_{t}^{(1/2)}\right)^{1/2}, S_{t}^{(1/2)}\le \delta-y]\P(D_{t}^{\phi_{1}}\in dy)\\
\geq& \E[\left( S_{t}^{(1/2)}\right)^{1/2}, S_{t}^{(1/2)}\le \delta/2] 
\times \P(D_{t}^{\phi_{1}}\le \delta/2).
\end{align*}
It follows from \cite[Lemma 3.2]{ParkSong19} 
and the stochastic continuity of the subordinator $(D_t^{\Phi_1})$ that
$$
\liminf_{t\downarrow 0}\left(\frac{\E[\left( S_{t}^{(1/2)}\right)^{1/2}, S_{t}^{(1/2)}\le\delta/2]}{t\log(1/t)}\times \P(D_{t}^{\phi_{1}}\le \delta/2)  \right)=\frac{1}{\sqrt{\pi}},
$$
which yields the desired lower bound.
\qed

\begin{thm}\label{thm:withlowerorder}
Let $\Omega$ be a bounded open interval when $d=1$, or a bounded connected $C^{1,1}$ open set when $d\geq 2$.
Let $(W_t)$ be a Brownian motion independent of a subordinator $(D_t)$ with Laplace exponent of the form $\phi(\lam)=\lam^{1/2} +\phi_{1}(\lam)$, 
where $\phi_1$ is the Laplace exponent of a subordinator whose L\'evy measure has an almost monotone density and
$\phi_{1} \in\mathcal{R}_{\beta_1}(\infty)$ with $\beta_1\in (0,\frac12)$.
Then 
$$
\lim_{t\downarrow 0}\frac{|\Omega|-\tilde{Q}^{W\circ D}_{\Omega}(t)}{t\log(1/t)}
=\frac{2|\partial \Omega|}{\pi}.
$$
\end{thm}

\pf
Note that the asymptotic relation in \eqref{0016} is still valid in the setting when $\phi(\lam)=\lam^{1/2}+\phi_1(\lam)$. Combining \eqref{0016} and Lemma \ref{lemma:mixedconv} yields conditions \eqref{0009} and \eqref{eqn:independent limit}, and the exact asymptotic rate of the spectral heat content follows from Proposition \ref{prop:general_timechange}.
\qed

\subsection{Case 3: $\phi\in \mathcal{R}_{\beta}(\infty)$ with $\beta\in (0,1/2)$}
Now we consider the case when $\phi\in \mathcal{R}_{\beta}(\infty)$ with $\beta\in (0,\frac12)$. First we establish the continuity of the function $t\mapsto \tilde{Q}^{W\circ D}_{\Omega}(t)$.

\begin{lemma}\label{lemma:continuous}
Let $\Omega$ be a bounded open interval when $d=1$, or a bounded connected $C^{1,1}$ open set when $d\geq 2$.
Let $(W_t)$ be a Brownian motion independent of a subordinator $(D_t)$ with infinite L\'evy measure.
Then the map $t\mapsto \tilde{Q}^{W\circ D}_{\Omega}(t)$ is continuous on $(0,\infty)$.
\end{lemma}
\pf
For any $t,h>0$,
\begin{eqnarray*}
0\leq \tilde{Q}^{W\circ D}_{\Omega}(t)-\tilde{Q}^{W\circ D}_{\Omega}(t+h)=\int_{\Omega}\P_{x}\left(D_{t}<\tau_{\Omega}^{BM}\leq D_{t+h}\right)dx.
\end{eqnarray*}
By the right-continuity of sample paths of the L\'evy process $(D_{t})$ and the continuity of the probability measure, 
$$
\lim_{h\downarrow 0}\P_{x}(D_{t}<\tau_{\Omega}^{BM}\leq D_{t+h})=\P_{x}(D_{t}< \tau_{\Omega}^{BM}\leq D_{t})
=\P_x(\emptyset)=0.
$$
Since $|\Omega|<\infty$, the dominated convergence theorem yields 
$$
\lim_{h\downarrow 0}\tilde{Q}^{W\circ D}_{\Omega}(t+h)=\tilde{Q}^{W\circ D}_{\Omega}(t).
$$

On the other hand,  note that for any $h>0$,
\begin{eqnarray*}
0\leq\tilde{Q}^{W\circ D}_{\Omega}(t-h)-\tilde{Q}^{W\circ D}_{\Omega}(t)=\int_{\Omega}\P_{x}\left(D_{t-h}<\tau_{\Omega}^{BM}\leq D_{t}\right)dx.
\end{eqnarray*}
It is known that a given L\'evy process $(D_{t})$ 
satisfies $\P^D(D_{t-}=D_{t})=1$ for any fixed $t>0$ (see e.g.\ \cite[Lemma 2.3.2]{Applebaum}). 
Thus,
$$
\lim_{h\downarrow 0}\P_{x}\left(D_{t-h}<\tau_{\Omega}^{BM}\leq D_{t}\right)
=\P_{x}(D_{t-}\leq \tau_{\Omega}^{BM}\leq D_{t})=\P_{x}(\tau_{\Omega}^{BM}=D_t)
=\E^D[\P^W_{x}(\tau_{\Omega}^{BM}=D_t)]. 
$$
The latter equals 0 since $\P^W_{x}(\tau_{\Omega}^{BM}=u)=0$ for all $x\in \R^{d}$ and $u>0$ due to \cite[Proposition 1.20]{Chung:1995}.
Since $|\Omega|<\infty$, the dominated convergence theorem yields
$$
\lim_{h\downarrow 0}\tilde{Q}^{W\circ D}_{\Omega}(t-h)=\tilde{Q}^{W\circ D}_{\Omega}(t),
$$
which completes the proof. 
\qed

The following proposition on weak convergence is needed to derive a result on the spectral heat content for the time-changed Brownian motion $(W_{D_t})$ when $\phi\in \mathcal{R}_\beta(\infty)$ with $\beta\in(0,\frac 12)$.
The proof is omitted since the statement follows easily by combining \cite[Lemma 1]{CG18} and \cite[Equation (11)]{KuhnSchilling2019} (see also \cite[Theorem 3.1]{KuhnSchilling2019}). 
Note that related statements appear in \cite[Theorem 8.7]{Sato} and \cite[Proposition 4.13]{Rosinski_isomorphism}.

\begin{prop}\label{prop:1/2}
Suppose $(D_{t})$ is a subordinator with Laplace exponent $\phi\in \mathcal{R}_{\beta}(\infty)$ with $\beta \in (0,1)$ and its L\'evy measure has an almost monotone density.
Let $f$ be a bounded continuous function on $(0,\infty)$ such that $\lim_{x\downarrow 0}\frac{f(x)}{x^{\gamma}}$ exists as a finite number for some constant $\gamma>\beta$. Then
$$
\lim_{t\downarrow 0}\int_{0}^{\infty}f(x)\frac{\P(D_{t}\in dx)}{t}=\int_{0}^{\infty}f(x)\nu(dx),
$$
where $\nu$ is the L\'evy measure of $(D_t)$.
\end{prop}
\begin{remark}
\begin{em}
A result similar to Proposition \ref{prop:1/2} is found in \cite[Proposition 3.4]{ParkSong19}, but the proof given in the latter paper involves a minor error. Indeed, since $t_{0}$ in that proof depends not only $\eta$ but also $\eps$, it is unclear whether it is possible to let $\eps\to 0$ uniformly for all $t\leq t_{0}$. To avoid this issue, one can instead use Proposition \ref{prop:1/2} above. 
\end{em}
\end{remark}

\begin{thm}\label{thm:lessthan1/2}
Let $\Omega$ be a bounded open interval when $d=1$, or a bounded connected $C^{1,1}$ open set when $d\geq 2$.
Suppose $(D_{t})$ is a subordinator with Laplace exponent $\phi\in \mathcal{R}_{\beta}(\infty)$ with $\beta \in (0,\frac 12)$ and its L\'evy measure has an almost monotone density.
Let $(W_t)$ be a Brownian motion independent of $(D_t)$.
Then
$$
\lim_{t\downarrow 0}\frac{|\Omega|-\tilde{Q}^{W\circ D}_{\Omega}(t)}{t}=\int_{0}^{\infty}(|\Omega|-Q_{\Omega}^{W}(u))\nu(du).
$$
\end{thm}

\pf
Let $f(x):=|\Omega|-Q_{\Omega}^{W}(x)$ for $x>0$. Then $f$ is bounded on $(0,\infty)$ since the domain $\Omega$ is bounded. Moreover, $f$ is continuous on $(0,\infty)$ due to Lemma \ref{lemma:continuous}. The restriction that $\beta\in(0,\frac 12)$ allows us to take $\gamma\in(\beta,\frac12)$, and for this $\gamma$, $\lim_{x\downarrow 0}f(x)/x^\gamma$ exists and equals 0 due to the asymptotic representation of $f$ in \eqref{0031}.  The desired result now follows upon noting that 
\[
	\frac{|\Omega|-\tilde{Q}^{W\circ D}_{\Omega}(t)}{t}
	=\frac{\E[f(D_t)]}{t}
	=\int_{0}^{\infty}f(x)\frac{\P(D_{t}\in dx)}{t}
\]
and applying Proposition \ref{prop:1/2}. 
\qed

\section{Time change by inverse subordinators}\label{section:inverse subordinators}
We now turn our attention to subdiffusion processes that are given as time-changed Brownian motions of the form $(W_{E_t})$, where the time change $(E_t)$  
is the (generalized) inverse of a subordinator $(D_t)$; i.e.,\ $E_t:=\inf\{u>0:D_u>t\}$. 
Here, we consider a general subordinator, assuming only that it has Laplace exponent $\phi\in \mathcal{R}_\beta(\infty)$ with $\beta\in(0,1)$; i.e.,\
\begin{align}\label{0017}
	\E[e^{-sD_t}]=e^{-t\phi(s)} \ \ \textrm{with} \ \ \phi(s)=s^\beta \ell(s)
\end{align}
for some function $\ell$ that is slowly varying at $\infty$ (i.e.,\ $\ell\in \mathcal{R}_0(\infty)$). 
Note that $(D_t)$ has strictly increasing, c\`adl\`ag paths with $D_0=0$, $0<D_t<\infty$ for all $t>0$, and $D_t\to \infty$ as $t\to\infty$, which implies that $(E_t)$ has nondecreasing, continuous paths with $E_0=0$ and $0<E_t<\infty$ for all $t>0$.

\subsection{Spectral and regular heat contents for Brownian motions time-changed by inverse subordinators}

Below we resort to Tauberian theorems of both Karamata's (polynomial) type  and de Bruijn's (exponential) type together with their associated monotone density theorems.  
The following lemma is a corollary to de Bruijn's Tauberian theorem in \cite[Theorem 4.12.9]{Bingham_book}; the proof is postponed to Appendix.

\begin{lemma}\label{lemma:de Brujin}
Let $\mu$ be a measure on $(0,\infty)$ whose Laplace transform $M(s):=\int_{0}^{\infty}e^{-s t}\mu(dt)$ converges for all $s>0$. Let $\beta\in(0,1)$. 
Then 
\beq\label{eqn:de Brujin1}
-\log M(s) \sim c_1 s^{\beta}\ell(s) \ \text{ as } \ s\to \infty 
\eeq
for some $\ell\in \mathcal{R}_0(\infty)$ and some constant $c_1>0$ if and only if
\beq\label{eqn:de Brujin2}
-\log \mu(0,t]\sim c_2 t^{-\frac{\beta}{1-\beta}}\ell_{1}(t) \ \text{ as } \ t\downarrow 0
\eeq
for some $\ell_1\in \mathcal{R}_0(0^+)$ and some constant $c_2>0$.
\end{lemma}

The following proposition is a key to establishing the main result of this section in Theorem \ref{thm:inverse_subordinator}. Note that the proof would be substantially simplified if we only considered the stable case; see Remark \ref{remark:inversestable} after the proof. 
\begin{prop}\label{prop:E_t}
Let $(E_t)$ be the inverse of a subordinator $(D_t)$ with Laplace exponent 
$\phi\in \mathcal{R}_\beta(\infty)$ with $\beta\in (0,1)$.
Then for any fixed $p>0$ and $\delta>0$, 
\[
	\E[E_t^p]\sim \E[E_t^p\mathbf{1}_{\{E_t\le \delta\}}]
	\sim \frac{\Gamma(p+1)}{\Gamma(p\beta+1)} 
	[\phi(t^{-1})]^{-p}
	 \ \ \textrm{as} \ \ t\downarrow 0.
\]
In particular, condition \eqref{eqn:independent limit} holds for the inverse subordinator $(E_t)$.
\end{prop}

\pf Express $\phi$ as $\phi(s)=s^\beta \ell(s)$ with $\ell\in \mathcal{R}_0(\infty)$, as in \eqref{0017}. 
We first derive 
 \begin{align}\label{0005}
 	\E[E_t^p]\sim \frac{\Gamma(p+1)}{\Gamma(p\beta+1)}t^{p\beta}[\ell(t^{-1})]^{-p} \ \ \textrm{as} \ \ t\downarrow 0. 
\end{align} 
By the Fubini Theorem and \eqref{0017}, 
\begin{align}
	&\int_0^\infty e^{-st}\P(D_x<t)dt dx
	=\int_{0}^{\infty}e^{-st}\int_{0}^{t}\P(D_{x}\in du)dt\notag\\
	=&\int_{0}^{\infty}\int_{u}^{\infty}e^{-st}dt\P(D_{x}\in du)
	=\frac{\E[e^{-sD_x}]}{s}=\frac{e^{-x\phi(s)}}{s}. \label{0018}
\end{align}
This, together with the formula $\E[X^p]=\int_0^\infty px^{p-1}\P(X>x)dx$ valid for any constant $p>0$ and any random variable $X\ge 0$, yields
\begin{align*}
&\int_0^\infty e^{-st}\E[E_t^p]dt
	=\int_0^\infty e^{-st}\left(\int_0^\infty p x^{p-1}\P(E_t>x)dx\right) dt\\
	=&\int_0^\infty px^{p-1} \left(\int_0^\infty e^{-st}\P(D_x<t)dt\right) dx
	=\frac ps\int_0^\infty x^{p-1} e^{-x\phi(s)} dx. 
\end{align*}
By the change of variables $y=x\phi(s)$, 
\begin{align*}
	\int_0^\infty e^{-st}\E[E_t^p]dt
	=\frac ps\int_0^\infty \left(\frac{y}{\phi(s)}\right)^{p-1} e^{-y} \frac{dy}{\phi(s)}
	=\frac{p\Gamma(p)}{s[\phi(s)]^p}
	=\Gamma(p+1)s^{-(p\beta+1)}[\ell(s)]^{-p}. 
\end{align*}
Application of Karamata's Tauberian theorem (see \cite[Theorem 1.7.1]{Bingham_book}) yields 
\[
	\int_0^t \E[E_r^p] dr \sim \frac{\Gamma(p+1)}{\Gamma(p\beta+2)}t^{p\beta+1}[\ell(t^{-1})]^{-p} \ \ \textrm{as} \ \ t\downarrow 0, 
\]
but since the function $t\mapsto \E[E_t^p]$ is nondecreasing, the monotone density theorem (see \cite[Theorem 1.7.2b]{Bingham_book}) applies and gives \eqref{0005}.

We now turn our attention to the asymptotic behavior of $t\mapsto\E[E_t^p\mathbf{1}_{\{E_t\le \delta\}}]$. Note that since the sample paths of the two processes $(E_t^p)_{t\ge 0}$ and $(\mathbf{1}_{\{E_t\le \delta\}})_{t\ge 0}$ are nondecreasing and nonincreasing, respectively, it is unclear whether or not the function $t\mapsto\E[E_t^p\mathbf{1}_{\{E_t\le \delta\}}]$ is monotone, and hence, an argument similar to the one provided above for the monotone function $t\mapsto\E[E_t^p]$ would not work. However, a similar discussion is still applicable to the function $t\mapsto\E[E_t^p\mathbf{1}_{\{E_t> \delta\}}]$, which is nondecreasing (notice the direction of the inequality). 
Indeed, we claim that the latter function possesses an exponential decay of the form
\begin{align}\label{0001}
	-\log \E[E_t^p\mathbf{1}_{\{E_t> \delta\}}] \sim c t^{-\beta/(1-\beta)}\ell_1(t)\ \ \textrm{as} \ \ t\downarrow 0
\end{align}
with some constant $c>0$ and some function $\ell_1\in \mathcal{R}_0(0^+)$. This asymptotic result together with \eqref{0005} implies
\[
	\lim_{t\downarrow 0}\frac{\E[E_t^p\mathbf{1}_{\{E_t\le  \delta\}}]}{t^{p\beta}[\ell(t^{-1})]^{-p}}
	=\lim_{t\downarrow 0}\frac{\E[E_t^p]}{t^{p\beta}[\ell(t^{-1})]^{-p}}
	- \lim_{t\downarrow 0}\frac{\E[E_t^p\mathbf{1}_{\{E_t> \delta\}}]}{t^{p\beta}[\ell(t^{-1})]^{-p}}
	=\frac{\Gamma(p+1)}{\Gamma(p\beta+1)},
\]
thereby completing the proof of the proposition.

Observe that
\begin{align*}
	\int_0^\infty e^{-st}\E[E_t^p\mathbf{1}_{\{E_t> \delta\}}]dt
	=\int_0^\infty e^{-st}\E[(E_t\mathbf{1}_{\{E_t> \delta\}})^p]dt
	=\int_0^\infty e^{-st}\int_0^\infty p x^{p-1}\P(E_t\mathbf{1}_{\{E_t> \delta\}}>x)dx dt.
\end{align*}
Note that $\{E_t\mathbf{1}_{\{E_t> \delta\}}>x\}=\{E_t>\max\{\delta,x\}\}$ for any fixed $x\ge 0$, so the latter integral can be decomposed as
\begin{align*}
	\int_0^\infty e^{-st}\int_0^\delta p x^{p-1}\P(E_t> \delta)dx dt
		+\int_0^\infty e^{-st}\int_\delta^\infty p x^{p-1}\P(E_t>x)dx dt
	=: J_1(s)+J_2(s).
\end{align*}
Now, by the Fubini theorem and \eqref{0018}, 
\begin{align*}
	J_1(s)
	&=\left(\int_0^\delta p x^{p-1} dx\right)\left(\int_0^\infty e^{-st}\P(D_\delta<t) dt\right)
	=\frac{\delta^p e^{-\delta\phi(s)}}{s},
\end{align*}
and
\begin{align*}
	J_2(s)
	&=\int_\delta^\infty  p x^{p-1} \left(\int_0^\infty e^{-st}\P(D_x<t) dt\right) dx
	=\int_\delta^\infty p x^{p-1} \frac{e^{-x\phi(s)}}{s}dx\\
	&=\frac ps\int_{\delta\phi(s)}^\infty \left(\frac{y}{\phi(s)}\right)^{p-1} e^{-y}\frac{dy}{\phi(s)}
	=\frac{p\Gamma(p,\delta\phi(s))}{s[\phi(s)]^p}, 
\end{align*}
where $\Gamma(a,z)$ denotes the upper incomplete Gamma function.
Putting the above together and using the identity $\Gamma(p+1,z)=p\Gamma(p,z)+z^p e^{-z}$ yields
\begin{align}
	\int_0^\infty e^{-st}\E[E_t^p\mathbf{1}_{\{E_t> \delta\}}]dt
	=\frac{\delta^p e^{-\delta\phi(s)}}{s}+\frac{p\Gamma(p,\delta\phi(s))}{s[\phi(s)]^p}
	= \frac{\Gamma(p+1,\delta\phi(s))}{s[\phi(s)]^p}\label{0002}. 
\end{align}
 Since $\Gamma(a,z)\sim z^{a-1}e^{-z}$ as $z\to\infty$ along the positive real line (see \cite[Formula 6.5.32]{Abramowitz}), it follows from \eqref{0002} that  
\[
	-\log\left(\int_0^\infty e^{-st}\E[E_t^p\mathbf{1}_{\{E_t> \delta\}}] dt\right) \sim \delta \phi(s)=\delta s^\beta\ell(s)\ \ \textrm{as} \ \ s\to\infty. 
\]
Thus, by Lemma \ref{lemma:de Brujin}, 
\[
	-\log\left(\int_0^t \E[E_r^p\mathbf{1}_{\{E_r> \delta\}}] dr\right) \sim c t^{-\beta/(1-\beta)}\ell_1(t) \ \ \textrm{as} \ \ t\downarrow 0
\]
for some constant $c>0$ and some function $\ell_1\in \mathcal{R}_0(0^+)$.
Setting $u=1/r$ and $x=1/t$ in the latter gives  
\[
	-\log\left(\int_x^\infty \frac{1}{u^2} \E[E_{1/u}^p\mathbf{1}_{\{E_{1/u}> \delta\}}] du\right)
	=-\log\left(\int_x^\infty e^{-h(u)} du\right) 
	\sim c x^{\beta/(1-\beta)}\ell_1(x^{-1}) \ \ \textrm{as} \ \ x\to \infty, 
\]
where
\[
	h(x):=-\log\left(\frac{1}{x^2} \E[E_{1/x}^p\mathbf{1}_{\{E_{1/x}> \delta\}}]\right). 
\]
Since $h(x)$ is nondecreasing, the monotone density theorem in \cite[Theorem 4.12.10(i)]{Bingham_book} is applicable and gives $h(x)\sim c x^{\beta/(1-\beta)}\ell_1(x^{-1})$ as $x\to\infty$. In other words,  
\[
	-\log\left(t^2 \E[E_t^p\mathbf{1}_{\{E_t> \delta\}}]\right)\sim ct^{-\beta/(1-\beta)}\ell_1(t) \ \ \textrm{as} \ \ t\downarrow 0,
\]
which in turn yields \eqref{0001}, as desired. 
\qed

\begin{remark}\label{remark:inversestable}
\begin{em}
Proposition \ref{prop:E_t} follows immediately if $(D_t)$ is a \textit{stable} subordinator $(S_t^{(\beta)})$ with index $\beta\in(0,1)$, in which case $\phi(s)=s^\beta$.
Indeed, the self-similarity of $(S_t^{(\beta)})$ with index $1/\beta$ implies the self-similarity of the inverse $(E_t^{(\beta)})$ with index $\beta$ (see e.g.\ \cite[Proposition 3.1]{MS_1}), and the $p$th moment of the inverse $\beta$-stable subordinator $(E_t^{(\beta)})$ is well-known and given by 
\begin{align}\label{0022}
	\E[(E^{(\beta)}_t)^p]=t^{p\beta} \E[(E^{(\beta)}_1)^p]=t^{p\beta}\frac{\Gamma(p+1)}{\Gamma(p\beta+1)} 
	\ \ \textrm{for all} \ \ t>0
\end{align}
(see e.g.\ \cite[Proposition 5.6]{HKU-book}).
Moreover, by the self-similarity again and the monotone convergence theorem, 
\[
	\E[(E^{(\beta)}_t)^p\mathbf{1}_{\{E^{(\beta)}_t\le \delta\}}]
	=t^{p\beta} \E[(E^{(\beta)}_1)^p\mathbf{1}_{\{E^{(\beta)}_1\le \delta t^{-\beta}\}}]
	\sim t^{p\beta} \E[(E^{(\beta)}_1)^p] \ \ \textrm{as} \ \ t\downarrow 0,
\]
which yields Proposition \ref{prop:E_t}. 
\end{em}
\end{remark}

Combining Proposition \ref{prop:general_timechange} with Proposition \ref{prop:E_t} yields the following main result of this section. 
Recall that the spectral heat content for a subordinated killed Brownian motion and that for a killed subordinated Brownian motion are defined in \eqref{0053}.

\begin{thm}\label{thm:inverse_subordinator}
Let $(W_t)$ be a Brownian motion independent of a subordinator $(D_t)$ with Laplace exponent 
$\phi\in \mathcal{R}_\beta(\infty)$ with $\beta\in (0,1)$.
Let $(E_t)$ be the inverse of $(D_t)$. 
Let $\Omega$ be a bounded open interval when $d=1$, or a bounded connected $C^{1,1}$ open set when $d\geq 2$.
Then the spectral heat content for the subordinated killed Brownian motion and that for the killed subordinated Brownian motion coincide; i.e., for any fixed $t>0$,
\begin{align}\label{0020}
	\tilde{Q}^{W\circ E}_{\Omega}(t)=Q^{W\circ E}_{\Omega}(t).
\end{align}
Moreover, 
\begin{align}\label{0019}
	\lim_{t\downarrow 0}\frac{|\Omega|-Q^{W\circ E}_{\Omega}(t)}
	{[\phi(t^{-1})]^{-1/2}}
	=\frac{|\partial \Omega|}{\Gamma(\beta/2+1)}
\ \  \ \textrm{and} \ \ \ 
	\lim_{t\downarrow 0}\frac{|\Omega|-\HH^{W\circ E}_{\Omega}(t)}{[\phi(t^{-1})]^{-1/2}}
	=\frac{|\partial \Omega|}{2\Gamma(\beta/2+1)}.
\end{align}
\end{thm}

\pf
Recall that $(E_t)$ has continuous paths. 
Since the notions of a killed subordinated Brownian motion and a subordinated killed Brownian motion coincide when the time change has continuous paths, the statement \eqref{0020} follows.

Express $\phi$ as $\phi(s)=s^\beta \ell(s)$ and $\ell\in \mathcal{R}_0(\infty)$, as in \eqref{0017}. 
For a fixed $\delta>0$, since $-\log \E[e^{-sD_\delta}]=\delta s^\beta\ell(s)$ for all $s>0$,
it follows from Lemma \ref{lemma:de Brujin} that 
\begin{align}\label{0006}
	-\log \P(D_\delta\le t)\sim ct^{-\beta/(1-\beta)}\ell_2(t) \ \ \textrm{as} \ \ t\downarrow 0 
\end{align}
for some constant $c>0$ and some function $\ell_2\in \mathcal{R}_0(0^+)$. In other words, the small ball probability $\P(D_\delta\le t)$ decays exponentially as $t\downarrow 0$. This together with Proposition \ref{prop:E_t} implies that 
\begin{align*}
	\lim_{t\downarrow 0}\frac{\P(E_t> \delta)}{\E[\sqrt{E_t}\mathbf{1}_{\{E_t\le \delta\}}]}
	=\lim_{t\downarrow 0}\frac{\P(D_\delta<t)}{\E[\sqrt{E_t}\mathbf{1}_{\{E_t\le \delta\}}]}
	=0.
\end{align*}
Thus, condition \eqref{0009} holds for the time change $(E_t)$. 
The desired result \eqref{0019} now follows by combining 
Propositions \ref{prop:general_timechange} and \ref{prop:E_t}.
\qed

\begin{remark}
\begin{em}
As in Remark \ref{remark:interpretation1}, we can re-express \eqref{0019} using a \textit{one-dimensional} Brownian motion $(B_t)$ independent of the inverse $(E^{(\beta)}_t)$ of a $\beta$-stable subordinator $(S_{t}^{(\beta)})$ as 
\begin{align}\label{0055}
	\lim_{t\downarrow 0}\frac{|\Omega|-Q^{W\circ E}_{\Omega}(t)}
	{[\phi(t^{-1})]^{-1/2}}
	=\E[\sup_{0\le u\leq E^{(\beta)}_{1}}B_{u}]|\partial \Omega|
\ \  \ \textrm{and} \ \ \ 
	\lim_{t\downarrow 0}\frac{|\Omega|-\HH^{W\circ E}_{\Omega}(t)}{[\phi(t^{-1})]^{-1/2}}
	=\frac12\E[\sup_{0\le u\leq E^{(\beta)}_{1}}B_{u}]|\partial \Omega|. 
\end{align}
Indeed, combining \eqref{0054}, \eqref{0022} with $p=1/2$, and the independence assumption gives
\[
	\E[\sup_{0\le u\leq E_{1}^{(\beta)}}B_{u}]
	=\E[(E_{1}^{(\beta)})^{1/2}]\times \frac{2}{\sqrt{\pi}}
	=\frac{1}{\Gamma(\beta/2+1)},
\]
and hence, \eqref{0019} and \eqref{0052} are equivalent. 
\end{em}
\end{remark}

\subsection{A complete asymptotic expansion}
In this short section, we derive a complete asymptotic expansion of the spectral heat content of a time-changed Brownian motion $(W_{E^{(\beta)}_t})$, where the time change $(E^{(\beta)}_t)$ is the inverse of a $\beta$-stable subordinator $(S^{(\beta)}_t)$. Recall that $E^{(\beta)}_t$ has moments of all orders given in \eqref{0022} and that the small ball probability of $S^{(\beta)}_\delta$ for a fixed $\delta>0$ decays exponentially as observed in \eqref{0006}.

It has been proved in \cite{vandenBergGilkey94} that, if $\Omega$ is a compact Riemannian manifold with a $C^{\infty}$ boundary,
then the spectral heat content $Q_{\Omega}^{W}(t)$ of a Brownian motion $(W_t)$ has the asymptotic expansion
\begin{equation*}
|\Omega|-Q_{\Omega}^{W}(t)\sim \sum_{n=1}^{\infty}c_{n}t^{\frac{n}{2}} \ \  \text{as} \ \    t\downarrow 0,
\end{equation*}
where $c_n$ are some constants.
The latter expression is understood as follows: for any fixed $N\in\mathbb{N}$, 
\begin{align}\label{0025}
|\Omega|-Q_{\Omega}^{W}(t)=\sum_{n=1}^{N}c_{n}t^{\frac{n}{2}} +O(t^{\frac{N+1}{2}})  \ \  \text{as} \ \    t\downarrow 0. 
\end{align}

\begin{thm}\label{thm:complete expansion}
Let $(E^{(\beta)}_t)$ be the inverse of a stable subordinator $(S^{(\beta)}_t)$ with index $\beta\in(0,1)$ which is independent of a Brownian motion $(W_t)$. 
Let $\Omega$ be a bounded connected open set in $\R^{d}$ with a $C^\infty$ boundary.
Then the spectral heat content of the time-changed Brownian motion $(W_{E^{(\beta)}_t})$ satisfies the asymptotic expansion
$$
|\Omega|-\tilde{Q}^{W\circ E^{(\beta)}}_{\Omega}(t)\sim \sum_{n=1}^{\infty}c_{n}\frac{\Gamma(1+\frac{n}{2})}{\Gamma(1+\frac{n\beta}{2})}t^{\frac{\beta n}{2}} \ \ \text{as} \ \  t\downarrow 0,
$$
which means that 
for any fixed $N\in \mathbb{N}$, 
\begin{align}\label{0024}
|\Omega|-\tilde{Q}^{W\circ E^{(\beta)}}_{\Omega}(t)=\sum_{n=1}^{N}c_{n}\frac{\Gamma(1+\frac{n}{2})}{\Gamma(1+\frac{n\beta}{2})}t^{\frac{\beta n}{2}} +O(t^{\frac{\beta(N+1)}{2}})  \ \ \text{as} \ \   t\downarrow 0.
\end{align}
\end{thm}

\pf
For simplicity of notation, we express $E^{(\beta)}_t$ as $E_t$ throughout the proof.
For a fixed $N\in \mathbb{N}$, by \eqref{0025}, there exist $\delta>0$ and $c>0$ such that 
\beq\label{eqn:complete aux1}
\left| |\Omega|-Q_{\Omega}^{W}(t) -\sum_{n=1}^{N}c_{n}t^{\frac{n}{2}}\right|\leq ct^{\frac{N+1}{2}}
\eeq
for all $t\in (0,\delta]$.
Hence, 
\begin{align}\label{eqn:complete aux2}
&|\Omega|-\tilde{Q}^{W\circ E}_{\Omega}(t)
	=\int_{\Omega}\P^W_{x}\times \P^{E}(\tau_{\Omega}^{BM}\leq E_{t})dx
	=\E^{E}\left[\int_{\Omega} \P^W_x(\tau_{\Omega}^{BM}\leq E_{t})dx \right]\nn\\
	=&\E^{E}[|\Omega|-Q_{\Omega}^{W}(E_t)]
	=\E^{E}[(|\Omega|-Q_{\Omega}^{W}(E_t))\mathbf{1}_{\{E_t\le \delta\}}]
		+\E^{E}[(|\Omega|-Q_{\Omega}^{W}(E_t))\mathbf{1}_{\{E_t> \delta\}}]. 
\end{align}
For the rest of the proof, we drop the superscript $E$ from the expectation sign $\E^{E}$.
As for the second term in \eqref{eqn:complete aux2}, by the exponential decay in \eqref{0006}, 
\begin{align*}
	\E[(|\Omega|-Q_{\Omega}^{W}(E_t))\mathbf{1}_{\{E_t> \delta\}}]
	\le |\Omega|\P(E_t> \delta)
	=|\Omega|\P(S^{(\beta)}_\delta<t)
	=o(t^{\frac{\beta(N+1)}{2}}).
\end{align*}
On the other hand, in terms of the first expression of \eqref{eqn:complete aux2}, 
if follows from \eqref{eqn:complete aux1} that
\begin{align}\label{eqn:complete aux3}
&\E[(|\Omega|-Q_{\Omega}^{W}(E_t))\mathbf{1}_{\{E_t\le \delta\}}]
\le \E\left[\left(\sum_{n=1}^{N}c_{n}E_t^{\frac{n}{2}}+cE_t^{\frac{N+1}{2}}\right)\mathbf{1}_{\{E_t\le \delta\}}\right]\nn\\
&\le \sum_{n=1}^{N}c_{n}\frac{\Gamma(1+\frac{n}{2})}{\Gamma(1+\frac{n\beta}{2})}t^{\frac{n\beta}{2}}
	-\E\left[\sum_{n=1}^{N}c_{n}E_t^{\frac{n}{2}}\mathbf{1}_{\{E_t> \delta\}}\right]
	+c\E\left[E_t^{\frac{N+1}{2}}\mathbf{1}_{\{E_t\le \delta\}}\right],
\end{align}
where we used \eqref{0022}.
As for the second term in \eqref{eqn:complete aux3}, by the Cauchy-Schwarz inequality, the moment formula \eqref{0022}, and the exponential decay in \eqref{0006}, 
$$
\E[E_t^{\frac{n}{2}}\mathbf{1}_{\{E_t> \delta\}}]
\leq \E[E_{t}^{n}]^{1/2}\P(E_{t}>\delta)^{1/2}
= \E[E_{t}^{n}]^{1/2}\P(S^{(\beta)}_\delta< t)^{1/2}
=o(t^{\frac{\beta(N+1)}{2}}).
$$
As for the third term in \eqref{eqn:complete aux3}, again by \eqref{0022},
\beq\label{eqn:complete aux5}
\E\left[E_t^{\frac{N+1}{2}}\mathbf{1}_{\{E_t\le \delta\}}\right]
\le \E[E_{t}^{\frac{N+1}{2}}]=O(t^{\frac{\beta(N+1)}{2}}).
\eeq
Combining \eqref{eqn:complete aux2}--\eqref{eqn:complete aux5} yields one direction of the equality in \eqref{0024}. The other direction is derived in a similar manner by considering a lower bound in \eqref{eqn:complete aux3} by means of \eqref{eqn:complete aux1}.
\qed

\section{Examples}\label{section:examples}

This section examines the asymptotic behavior of the spectral heat content for subordinate killed Brownian motions with various underlying subordinators. 
Throughout the section, we assume that $\Omega$ is a bounded open interval in $\R^{1}$ or a  bounded connected $C^{1,1}$ open set in $\R^{d}$ with $d\geq 2$. 
\begin{enumerate}[(1)]
\item Stable subordinators

The most typical example is a $\beta$-stable subordinator $(S_t^{(\beta)})$ with $\beta\in (0,1)$ whose Laplace exponent is $\phi(\lam)=\lam^{\beta}$. 
In this case, by Theorem \ref{thm:stable2}, Proposition \ref{prop:stable2} or Theorem \ref{thm:withlowerorder}, and Theorem \ref{thm:lessthan1/2},
$$
\begin{cases}
\displaystyle\lim_{t\downarrow 0}\frac{|\Omega|-\tilde{Q}_{\Omega}^{W\circ S^{(\beta)}}(t)}{t^{1/2\beta}}=\E[(S_{1}^{(\beta)})^{1/2}]\times \frac{2|\partial \Omega|}{\sqrt{\pi}}, & \beta\in (\frac12,1),\\
\displaystyle\lim_{t\downarrow 0}\frac{|\Omega|-\tilde{Q}_{\Omega}^{W\circ S^{(\beta)}}(t)}{t\log(1/t)}=\frac{2|\partial \Omega|}{\pi}, &\beta=\frac12,\\
\displaystyle\lim_{t\downarrow 0}\frac{|\Omega|-\tilde{Q}_{\Omega}^{W\circ S^{(\beta)}}(t)}{t}=\int_{0}^{\infty}(|\Omega|-Q_{\Omega}^{W}(u))\frac{\beta}{\Gamma(1-\beta)}u^{-1-\beta}du, &\beta\in (0,\frac12).
\end{cases}
$$

In terms of the regular heat content, define the perimeter $\text{Per}_{X}(\Omega)$ for a given L\'evy process $(X_t)$ with L\'{e}vy measure $\nu$ by  
$$
\text{Per}_{X}(\Omega)= \int_{\Omega}\int_{\Omega ^c-x}\nu(dy) dx.
$$ 
Applying Theorem \ref{thm:stable2} when $\beta\in(\frac12,1)$, Proposition \ref{prop:stable2} or Theorem \ref{thm:withlowerorder} when $\beta=\frac12$, and \cite[Theorem 3]{CG17} or \cite[Theorem 3.2]{GPS19} when $\beta\in (0,\frac12)$, we obtain the asymptotic behavior of the regular heat content as follows: 
$$
\begin{cases}
\displaystyle\lim_{t\downarrow 0}\frac{|\Omega|-\HH_{\Omega}^{W\circ S^{(\beta)}}(t)}{t^{1/2\beta}}
=\E[(S_{1}^{(\beta)})^{1/2}]\times \frac{|\partial \Omega|}{\sqrt{\pi}}, & \beta\in (\frac12,1),\\
\displaystyle\lim_{t\downarrow 0}\frac{|\Omega|-\HH_{\Omega}^{W\circ S^{(\beta)}}(t)}{t\log(1/t)}=\frac{|\partial \Omega|}{\pi}, &\beta=\frac12,\\
\displaystyle\lim_{t\downarrow 0}\frac{|\Omega|-\HH_{\Omega}^{W\circ S^{(\beta)}}(t)}{t}=\text{Per}_{W\circ S^{(\beta)}}(\Omega)=\int_{\Omega}\int_{\Omega ^c}
\frac{c(d,2\beta) dy dx}{|x-y|^{d+2\beta}}, &\beta\in (0,\frac12),
\end{cases}
$$
where $c(d, \alpha):= \frac{ \alpha  \,
\Gamma(\frac{d+\alpha}2)} {2^{1-\alpha} \,
\pi^{d/2}\Gamma(1-\frac{\alpha}2)}$. 
On the other hand, for the inverse stable subordinator $(E_{t}^{\beta})$, it follows from Theorem \ref{thm:inverse_subordinator} that 
$$
\lim_{t\downarrow 0}\frac{|\Omega|-\tilde{Q}_{\Omega}^{W\circ E^{(\beta)}}}{t^{\beta/2}}=\frac{|\partial \Omega|}{\Gamma(\beta/2+1)}, \quad \beta\in (0,1), 
$$
and 
$$
\lim_{t\downarrow 0}\frac{|\Omega|-\HH_{\Omega}^{W\circ E^{(\beta)}}}{t^{\beta/2}}=\frac{|\partial \Omega|}{2\Gamma(\beta/2+1)}, \quad \beta\in (0,1).
$$
One interesting observation is that the rate function for $|\Omega|-\tilde{Q}_{\Omega}^{W\circ S^{(\beta)}}(t)$ defined by 
$$
R_{\beta}(t)=
\begin{cases}
t^{1/2\beta} &\text{if } \beta\in (\frac12,1),\\
t\log(1/t) &\text{if } \beta=\frac12,\\
t &\text{if } \beta\in (0,\frac12),
\end{cases}
$$
satisfies
$$
R_{\beta}(t) =o(t^{\beta/2}) \ \text{ as } \ t\downarrow 0 \ \text{ for any } \  \beta\in(0,1).
$$
Thus, \textit{regardless of the value of $\beta\in(0,1)$, the rate of decay of the spectral heat content with the time change being the inverse stable subordinator $(E^{(\beta)}_t)$ is greater than that with the time change given by the stable subordinator $(S^{(\beta)}_t)$.}

\item Mixed stable subordinators

A mixed stable subordinator (or a mixture of stable subordinators) $(D^{MSS}_t)$ is defined to be the sum of independent stable subordinators with different indices, the Laplace exponent of which takes the form $\phi^{MSS}(\lam)=\sum_{i=1}^{n}\lam^{\beta_{i}}$ with $0<\beta_{1}<\beta_{2}<\cdots < \beta_{n}<1$. The density of the associated L\'evy measure is given by 
$$
\nu^{MSS}(u)=\sum_{i=1}^{n}\frac{\beta_{i}}{\Gamma(1-\beta_{i})}u^{-1-\beta_{i}}.
$$
Since $\phi^{MSS}(\lam)\sim \lam^{\beta_n}$ as $\lam\to\infty$, it follows that $(\phi^{MSS})^{-1}(s)\sim s^{1/\beta_n}$ as $s\to\infty$, and hence, 
$$
(\phi^{MMS})^{-1}(1/t)\sim t^{-1/\beta_{n}} \ \text{ as } \ t\downarrow 0.
$$
Thus, by Theorems \ref{thm:stable2}, \ref{thm:withlowerorder} and \ref{thm:lessthan1/2} as well as \cite[Theorem 3]{CG17} or \cite[Theorem 3.2]{GPS19}, the largest index $\beta_{n}$ determines the asymptotic behaviors of the spectral and regular heat contents as follows:
$$
\begin{cases}
\displaystyle\lim_{t\downarrow 0}\frac{|\Omega|-\tilde{Q}_{\Omega}^{W\circ D^{MSS}}(t)}{t^{1/2\beta_{n}}}=\E[(S_{1}^{(\beta_n)})^{1/2}]\times \frac{2|\partial \Omega|}{\sqrt{\pi}}, & \beta_{n}\in (\frac12,1),\\
\displaystyle\lim_{t\downarrow 0}\frac{|\Omega|-\tilde{Q}_{\Omega}^{W\circ D^{MSS}}(t)}{t\log(1/t)}=\frac{2|\partial \Omega|}{\pi}, &\beta_{n}=\frac12,\\
\displaystyle\lim_{t\downarrow 0}\frac{|\Omega|-\tilde{Q}_{\Omega}^{W\circ D^{MSS}}(t)}{t}=\int_{0}^{\infty}(|\Omega|-Q_{\Omega}^{W}(u))\nu^{MMS}(du),&\beta_{n}\in (0,\frac12),
\end{cases}
$$
and 
$$
\begin{cases}
\displaystyle\lim_{t\downarrow 0}\frac{|\Omega|-\HH_{\Omega}^{W\circ D^{MSS}}(t)}{t^{1/2\beta_{n}}}=\E[(S_{1}^{(\beta_n)})^{1/2}]\times \frac{|\partial \Omega|}{\sqrt{\pi}}, & \beta_n\in (\frac12,1),\\
\displaystyle\lim_{t\downarrow 0}\frac{|\Omega|-\HH_{\Omega}^{W\circ D^{MSS}}(t)}{t\log(1/t)}=\frac{|\partial \Omega|}{\pi}, &\beta_n=\frac12,\\
\displaystyle\lim_{t\downarrow 0}\frac{|\Omega|-\HH_{\Omega}^{W\circ D^{MSS}}(t)}{t}=\text{Per}_{W\circ D^{MSS}}(\Omega),
&\beta_n\in (0,\frac12).
\end{cases}
$$
On the other hand, by Theorem \ref{thm:inverse_subordinator}, the inverse $(E^{MSS}_{t})$ of the mixed stable subordinator $(D^{MSS}_t)$ satisfies    
$$
\lim_{t\downarrow 0}\frac{|\Omega|-\tilde{Q}_{\Omega}^{W\circ E^{MSS}}}{t^{\beta_{n}/2}}=\frac{|\partial \Omega|}{\Gamma(\beta_n/2+1)}, \quad \beta_n\in (0,1),
$$
and 
$$
\lim_{t\downarrow 0}\frac{|\Omega|-\HH_{\Omega}^{W\circ E^{MSS}}}{t^{\beta_n/2}}=\frac{|\partial \Omega|}{2\Gamma(\beta_n/2+1)}, \quad \beta_n\in (0,1).
$$

\item Tempered stable subordinators 

A subordinator $(D^{TSS}_t)$ with Laplace exponent $\phi^{TSS}(\lam)=(\lam+\theta)^{\beta}-\theta^\beta$ with $\beta\in (0,1)$ and $\theta>0$ is called an exponentially tempered (or tilted) stable subordinator. It behaves like a stable process in short time scale and a Gaussian process in large time scale. 
The Laplace exponent satisfies $\phi^{TSS}\in \mathcal{R}_\beta(\infty)\cap \mathcal{R}_1(0^+)$ and 
the density of the associated L\'evy measure is given by
\[
	\nu^{TSS}(u)=\frac{\beta}{\Gamma(1-\beta)}e^{-\theta u} u^{-1-\beta},
\]
 which explains how the parameter $\theta$ helps reduce (or temper) the jump sizes of a $\beta$-stable subordinator $(S^{(\beta)}_t)$. 
By Theorems \ref{thm:stable2},  \ref{thm:withlowerorder} and \ref{thm:lessthan1/2}  as well as \cite[Theorem 3]{CG17} or \cite[Theorem 3.2]{GPS19},
$$
\begin{cases}
\displaystyle\lim_{t\downarrow 0}\frac{|\Omega|-\tilde{Q}_{\Omega}^{W\circ D^{TSS}}(t)}{t^{1/2\beta}}=\E[(S_{1}^{(\beta)})^{1/2}]\times \frac{2|\partial \Omega|}{\sqrt{\pi}}, & \beta\in (\frac12,1),\\
\displaystyle\lim_{t\downarrow 0}\frac{|\Omega|-\tilde{Q}_{\Omega}^{W\circ D^{TSS}}(t)}{t\log(1/t)}=\frac{2|\partial \Omega|}{\pi}, &\beta=\frac12,\\
\displaystyle\lim_{t\downarrow 0}\frac{|\Omega|-\tilde{Q}_{\Omega}^{W\circ D^{TSS}}(t)}{t}=\int_{0}^{\infty}(|\Omega|-Q_{\Omega}^{W}(u))\nu^{TSS}(du),&\beta\in (0,\frac12),
\end{cases}
$$
and 
$$
\begin{cases}
\displaystyle\lim_{t\downarrow 0}\frac{|\Omega|-\HH_{\Omega}^{W\circ D^{TSS}}(t)}{t^{1/2\beta}}=\E[(S_{1}^{(\beta)})^{1/2}]\times \frac{|\partial \Omega|}{\sqrt{\pi}}, & \beta\in (\frac12,1),\\
\displaystyle\lim_{t\downarrow 0}\frac{|\Omega|-\HH_{\Omega}^{W\circ D^{TSS}}(t)}{t\log(1/t)}=\frac{|\partial \Omega|}{\pi}, &\beta=\frac12,\\
\displaystyle\lim_{t\downarrow 0}\frac{|\Omega|-\HH_{\Omega}^{W\circ D^{TSS}}(t)}{t}=\text{Per}_{W\circ D^{TSS}}(\Omega),
&\beta\in (0,\frac12).
\end{cases}
$$
Notice that, \textit{whether it is the spectral or regular heat content, the tempering factor $\theta>0$ appears in the asymptotic limit only in the case when $\beta\in(0,\frac 12)$} (through $\nu^{TSS}$ or $\text{Per}_{W\circ D^{TSS}}(\Omega)$). 
On the other hand, by Theorem \ref{thm:inverse_subordinator}, the inverse $(E^{TSS}_{t})$ of the exponentially tempered stable subordinator $(D^{TSS}_t)$ satisfies    
$$
\lim_{t\downarrow 0}\frac{|\Omega|-\tilde{Q}_{\Omega}^{W\circ E^{TSS}}}{t^{\beta/2}}=\frac{|\partial \Omega|}{\Gamma(\beta/2+1)}, \quad \beta\in (0,1),
$$
and 
$$
\lim_{t\downarrow 0}\frac{|\Omega|-\HH_{\Omega}^{W\circ E^{TSS}}}{t^{\beta/2}}=\frac{|\partial \Omega|}{2\Gamma(\beta/2+1)}, \quad \beta\in (0,1).
$$

\end{enumerate}

\section*{Appendix} 

\textbf{Proof of Lemma \ref{lemma:de Brujin}.} 
We apply \cite[Theorem 4.12.9]{Bingham_book} using the notations appearing in that theorem, so in particular, the function $\phi$ below does not represent a Laplace exponent. 

Note first that for any $f(t)\in \mathcal{R}_{-\gamma}(0^+)$ with $\gamma>0$ and $\tilde{f}(s):=f(s^{-1})\in \mathcal{R}_\gamma(\infty)$, 
\begin{align}\label{0021}
	f^{\leftarrow}(s)=\sup\{t:f(t)>s\}
	=\sup\left\{\frac 1u:\tilde{f}(u)>s\right\}
	=\frac{1}{\inf\{u:\tilde{f}(u)>s\}}
	=\frac{1}{\tilde{f}^{\leftarrow}(s)}.
\end{align}
Let $\ell\in \mathcal{R}_0(\infty)$ and define 
\[
	\psi(t):=t^{-1/\beta} \ell_{2}^\#(t^{-1/\beta}) \ \ \textrm{and} \ \ \phi(t):=t\psi(t),  
\]
where $\ell_{2}(s):=[\ell(s)]^{1/\beta}$ 
and $\ell_{2}^\#(s)$ is the de Bruijn conjugate of $\ell_2(s)$. Then it follows that $\psi\in \mathcal{R}_{-1/\beta}(0^+)$ and $\phi\in \mathcal{R}_{-(1-\beta)/\beta}(0^+)$. 
By the choice of $\psi(t)$,  
\[
	\tilde{\psi}(s):=\psi(s^{-1})=s^{1/\beta} \ell_{2}^\#(s^{1/\beta}),
\]
where $\ell_{2}^\#(s^{1/\beta})\in \mathcal{R}_0(\infty)$, so by \eqref{0021} and \cite[Proposition 1.5.15]{Bingham_book}, 
\[
	\frac{1}{\psi^\leftarrow(s)}
	=\tilde{\psi}^{\leftarrow}(s)
	\sim s^\beta [\ell_{2}(s)]^\beta=s^\beta \ell(s) \ \ \textrm{as} \ \ s\to\infty.
\]
Hence, condition \eqref{eqn:de Brujin1} takes the form 
$$
-\log M(s) \sim \frac{c_1}{\psi^{\leftarrow}(s)} \ \ \textrm{as} \ \ s\to\infty
$$ 
with $\psi\in \mathcal{R}_{-1/\beta}(0^+)$. 
Applying de Bruijn's Tauberian theorem in \cite[Theorem 4.12.9]{Bingham_book} with $\alpha:=-(1-\beta)/\beta<0$, we observe that the latter condition is equivalent to the condition that 
$$
-\log \mu(0,t] \sim \frac{c_{2}}{\phi^{\leftarrow}(1/t)} \ \text{ as } \ t\downarrow 0
$$
for some constant $c_{2}>0$. Thus, the proof will complete once we show that the latter condition takes the form in \eqref{eqn:de Brujin2}. 

To find an expression for $\phi^{\leftarrow}(1/t)$, 
note that
\[
	\tilde{\phi}(s)
	:=\phi(s^{-1})
	=s^{-1}\psi(s^{-1})
	=s^{-1+1/\beta} \ell_{2}^\#(s^{1/\beta})
	=s^{(1-\beta)/\beta} [\ell_{3}(s)]^{(1-\beta)/\beta},
\]
where $\ell_{3}(s):= [\ell_{2}^\#(s^{1/\beta})]^{\beta/(1-\beta)}$.
By \cite[Proposition 1.5.15]{Bingham_book} again, 
\[
	\tilde{\phi}^{\leftarrow}(s)\sim s^{\beta/(1-\beta)} \ell_{3}^\# (s^{\beta/(1-\beta)}), 
\]
so by \eqref{0021},
\[
	\frac{1}{\phi^{\leftarrow}(t^{-1})}
	=\tilde{\phi}^{\leftarrow}(t^{-1})
	\sim t^{-\beta/(1-\beta)} \ell_{3}^\# (t^{-\beta/(1-\beta)}) \  \ \textrm{as} \ \ t\downarrow 0.
\]
Letting $\ell_{1}(t):= \ell_{3}^\# (t^{-\beta/(1-\beta)}) \in \mathcal{R}_0(0^+)$ gives the form in \eqref{eqn:de Brujin2}, as desired. 
\qed

\noindent
\textbf{Acknowledgments:}
The authors are grateful to Professor Renming Song (University of Illinois at Urbana--Champign) for his valuable input regarding the material presented in Section 3. 
Kei Kobayashi's research was partially supported by a Faculty Fellowship at Fordham University.  
\bigskip

\bibliographystyle{plain}
\bibliography{KobayashiK}

\begin{thebibliography}{10}

\bibitem{Abramowitz}
M.~Abramowitz and I.~A. Stegun.
\newblock {\em Handbook of Mathematical Functions with Formulas, Graphs, and
  Mathematical Tables}.
\newblock Dover, 10th edition, 1972.

\bibitem{Valverde2014}
L.~Acu\~{n}a Valverde.
\newblock Trace asymptotics for fractional {S}chr\"{o}dinger operators.
\newblock {\em J. Funct. Anal.}, 266(2):514--559, 2014.

\bibitem{Valverde2017}
L.~Acu\~{n}a Valverde.
\newblock Heat content for stable processes in domains of {$\Bbb{R}^d$}.
\newblock {\em J. Geom. Anal.}, 27(1):492--524, 2017.

\bibitem{Applebaum}
D.~Applebaum.
\newblock {\em \textit{L{\'e}vy Processes and Stochastic Calculus}}.
\newblock Cambridge University Press, second edition, 2009.

\bibitem{Bertoin_Levy}
J.~Bertoin.
\newblock {\em L\'{e}vy {P}rocesses}.
\newblock Cambridge Tracts in Mathematics. Cambridge University Press, 1998.

\bibitem{Bingham_book}
N.~H. Bingham, C.~M. Goldie, and J.~L. Teugels.
\newblock {\em \textit{Regular Variation}}.
\newblock Encyclopedia of Mathematics and its Applications. Cambridge
  University Press, 1987.

\bibitem{Chung:1995}
K.~L. Chung and Z.~X. Zhao.
\newblock {\em From {B}rownian motion to {S}chr\"{o}dinger's equation}, volume
  312 of {\em Grundlehren der Mathematischen Wissenschaften [Fundamental
  Principles of Mathematical Sciences]}.
\newblock Springer-Verlag, Berlin, 1995.

\bibitem{CG17}
W.~Cygan and T.~Grzywny.
\newblock Heat content for convolution semigroups.
\newblock {\em J. Math. Anal. Appl.}, 446(2):1393--1414, 2017.

\bibitem{CG18}
W.~Cygan and T.~Grzywny.
\newblock A note on the generalized heat content for l\'evy processes.
\newblock {\em \textit{Bull. Korean Math. Soc.}}, 55(5):1463--1481, 2018.

\bibitem{FellerTwo71}
W.~Feller.
\newblock {\em An introduction to probability theory and its applications.
  {V}ol. {II}}.
\newblock Second edition. John Wiley \& Sons, Inc., New York-London-Sydney,
  1971.

\bibitem{Grzywny2018}
T.~Grzywny, L.~Le{\.z}aj, and B.~Trojan.
\newblock Transition densities of subordinators of positive order.
\newblock {\em https://arxiv.org/abs/1812.06793}, 2020.

\bibitem{GPS19}
T.~Grzywny, H.~Park, and R.~Song.
\newblock Spectral heat content for {L}\'{e}vy processes.
\newblock {\em Math. Nachr.}, 292(4):805--825, 2019.

\bibitem{JinKobayashi}
S.~Jin and K.~Kobayashi.
\newblock Strong approximation of stochastic differential equations driven by a
  time-changed {B}rownian motion with time-space-dependent coefficients.
\newblock {\em \textit{J. Math. Anal. Appl.}}, 476(2):619--636, 2019.

\bibitem{KS15}
K.~Kaleta and P.~Sztonyk.
\newblock Estimates of transition densities and their derivatives for jump
  {L}\'{e}vy processes.
\newblock {\em J. Math. Anal. Appl.}, 431(1):260--282, 2015.

\bibitem{KaratzasShreve}
I.~Karatzas and S.~Shreve.
\newblock {\em \textit{Brownian Motion and Stochastic Calculus}}.
\newblock Springer-Verlag New York, second edition, 1998.

\bibitem{Kobayashi_smallball}
K.~Kobayashi.
\newblock Small ball probabilities for a class of time-changed self-similar
  processes.
\newblock {\em \textit{Statist. Probab. Lett.}}, 110:155--161, 2016.

\bibitem{KuhnSchilling2019}
F.~K\"uhn and R.~Schilling.
\newblock On the domain of fractional {L}aplacians and related generators of
  {F}eller processes.
\newblock {\em \textit{J. Funct. Anal.}}, 276(8):2397--2439, 2019.

\bibitem{Magdziarz_spa}
M.~Magdziarz.
\newblock Stochastic representation of subdiffusion processes with
  time-dependent drift.
\newblock {\em \textit{Stoch. Proc. Appl.}}, 119:3238--3252, 2009.

\bibitem{MagdziarzZorawik}
M.~Magdziarz and T.~Zorawik.
\newblock {S}tochastic representation of fractional subdiffusion equation.
  {T}he case of infinitely divisible waiting times, {L}\'evy noise and
  space-time-dependent coefficients.
\newblock {\em \textit{Proc. Amer. Math. Soc.}}, 144:1767--1778, 2016.

\bibitem{MNV1}
M.~M. Meerschaert, E.~Nane, and P.~Vellaisamy.
\newblock {F}ractional {C}auchy problems on bounded domains.
\newblock {\em \textit{Ann. Probab.}}, 37:979--1007, 2009.

\bibitem{MS_1}
M.~M. Meerschaert and H-P. Scheffler.
\newblock Limit theorems for continuous-time random walks with infinite mean
  waiting times.
\newblock {\em \textit{J. Appl. Probab.}}, 41:623--638, 2004.

\bibitem{Miranda2007}
M.~Miranda~Jr, D.~Pallara, F.~Paronetto, and M.~Preunkert.
\newblock Short-time heat flow and fluctuations of bounded variation in
  $\mathbb{R}^{N}$.
\newblock {\em \textit{Ann. Fac. Sci. Toulouse}}, 16(1):125--145, 2007.

\bibitem{NaneNi2018}
E.~Nane and Y.~Ni.
\newblock Path stability of stochastic differential equations driven by
  time-changed {L}\'evy noises.
\newblock {\em \textit{ALEA, Lat. Am. J. Probab. Math. Stat.}}, 15:479--507,
  2018.

\bibitem{ParkSong19}
H.~Park and R.~Song.
\newblock Small time asymptotics of spectral heat contents for subordinate
  killed {B}rownian motions related to isotropic $\alpha$-stable processes.
\newblock {\em \textit{Bull. London Math. Soc.}}, 51:371--384, 2019.

\bibitem{ParkSong21}
H.~Park and R.~Song.
\newblock Spectral heat content $\alpha$-stable processes in $c^{1,1}$ open
  sets.
\newblock {\em Submitted. https://arxiv.org/abs/2007.02815}, 2021.

\bibitem{Rosinski_isomorphism}
J.~Rosi{\'n}ski.
\newblock Representations and isomorphism identities for infinitely divisible
  processes.
\newblock {\em \textit{Ann. Probab.}}, 46(6):3229--3274, 2018.

\bibitem{Sato}
K.~Sato.
\newblock {\em \textit{L{\'e}vy Processes and Infinitely Divisible
  Distributions}}.
\newblock Cambridge University Press, 1999.

\bibitem{Schilling}
R.~L. Schilling, R.~Song, and Z.~Vondracek.
\newblock {\em \textit{Bernstein Functions: Theory and Applications}}.
\newblock De Gruyter, 2010.

\bibitem{HKU-book}
S.~Umarov, M.~Hahn, and K.~Kobayashi.
\newblock {\em \textit{Beyond the {T}riangle: {B}rownian {M}otion, {I}t\^o
  {C}alculus, and {F}okker--{P}lanck {E}quation --- {F}ractional
  {G}eneralizations}}.
\newblock World Scientific, 2018.

\bibitem{vandenBerg1989}
M.~van~den Berg and E.~B. Davies.
\newblock Heat flow out of regions in {${\bf R}^m$}.
\newblock {\em Math. Z.}, 202(4):463--482, 1989.

\bibitem{vandenBergGilkey94}
M.~van~den Berg and P.~B. Gilkey.
\newblock Heat content asymptotics for a {R}iemannian manifold with boundary.
\newblock {\em \textit{J. Funct. Anal.}}, 120:48--71, 1994.

\bibitem{vandenBerg2015}
M.~van~den Berg and K.~Gittins.
\newblock Uniform bounds for the heat content of open sets in {E}uclidean
  space.
\newblock {\em Differential Geom. Appl.}, 40:67--85, 2015.

\bibitem{BergLeGall94}
M.~van~den Berg and J.-F. Le~Gall.
\newblock Mean curvature and the heat equation.
\newblock {\em Math. Z.}, 215(3):437--464, 1994.

\end{thebibliography}

\end{doublespace}

\vskip 0.3truein

{\bf Kei Kobayashi}

Department of Mathematics, Fordham University, NY 10023, USA

E-mail: \texttt{kkobayashi5@fordham.edu}

\vskip 0.3truein

{\bf Hyunchul Park}

Department of Mathematics, State University of New York at New Paltz, NY 12561,
USA

E-mail: \texttt{parkh@newpaltz.edu}

\bigskip

\end{document}